\documentclass[a4paper,10pt]{article}
\usepackage[margin=3cm]{geometry}
\usepackage{amsmath}
\usepackage{amsfonts}
\usepackage{amssymb}
\usepackage{graphicx}
\usepackage{verbatim}
\usepackage{changepage}
\usepackage{mathtools}
\usepackage{caption}
\usepackage{subcaption}
\usepackage{algorithm}
\usepackage{slashbox,hyperref}
\usepackage[square,numbers]{natbib}
\providecommand{\keywords}[1]
{
  \small	
  \textbf{\textit{Keywords---}} #1
}

\title{Multigrid preconditioning of singularly perturbed
  convection-diffusion equations}
\author{ M. Shahid$^{\nmid}$, S.P. MacLachlan$^{\ddagger}$ and H. bin Zubair$^{\nmid}$ \\
        \small $^{\nmid}$IBA-Karachi, Pakistan \\
        \small $^{\ddagger}$Memorial University of Newfoundland, Canada  \\
}
\date{} 

\newtheorem{example}{Example}[section]

\newcommand{\vb}{\boldsymbol{b}}
\newcommand{\vn}{\boldsymbol{n}}

\begin{document}
\maketitle
\begin{abstract}
Boundary value problems based on the convection-diffusion equation arise naturally in models of fluid flow across a variety of engineering applications and design feasibility studies. Naturally, their efficient numerical solution has continued to be an interesting and active topic of research for decades. In the context of finite-element discretization of these boundary value problems, the Streamline Upwind Petrov-Galerkin (SUPG) technique yields accurate discretization in the singularly perturbed regime. In this paper, we propose efficient multigrid iterative solution methods for the resulting linear systems. In particular, we show that techniques from standard multigrid for anisotropic problems can be adapted to these discretizations on both tensor-product as well as semi-structured meshes. The resulting methods are demonstrated to be robust preconditioners for several standard flow benchmarks.

\end{abstract}


\section{Introduction}\label{sec1}
 Convection-diffusion problems appear in many models related to fluid flow, including the Oseen equations, fluid transport, and related physical phenomena. They also appear when modeling studies based on the Navier-Stokes equations \cite[Chapter 1]{Morton96} are simplified under certain conditions. In the computational science and engineering communities, there has been a keen and unabated interest in their numerical solution for a considerably long time. This includes efforts aimed at developing accurate discretization schemes for these problems, as well as devising efficient numerical methods for the solution of the resulting linear systems. We are particularly interested in convection dominated models (i.e., those that reflect high P\'eclet number). This is the regime of problems for which standard numerical methods for elliptic equations are known to lead to inaccurate solutions that fail to resolve the boundary layers in the continuum solution and are prone to non-physical oscillations (``wiggles'') in the resulting numerical solutions~\cite{stynes2018convection}.  For this regime of problems, the use of properly chosen (``stable'') discretization schemes is essential to achieving accurate approximations; however, such schemes often complicate an iterative approach to the numerical solution. While there have been several recent papers on the stable upwind finite-difference discretization of such singularly perturbed convection-diffusion equations~\cite{echeverria2018convergence, etna_vol54_pp31-50, SMacLachlan_etal_2021a}, considerably less work appears in case of finite-element discretizations. In this paper, we develop efficient linear solvers for the well-known Streamline Upwind Petrov-Galerkin (SUPG) finite-element discretization \cite{hughes1987recent, hughes1987new} for singularly perturbed convection-diffusion equations discretized on layer-adapted meshes.
 
Consider a domain $\Omega \subset \mathbb{R}^2$ and the equation,
\begin{equation*}
- \varepsilon \Delta u - \vb(x,y) \cdot \nabla u + c(x,y)u = f,
\end{equation*}
posed on $\Omega$, with $\varepsilon \in (0,1]$, $c(x,y)\geq 0$, and $\vb \neq \boldsymbol{0}$.  The function, $u$, may represent the steady-state density of particles under weak diffusive effects being transported along by the flow velocity, $\vb$.  In the regime where the P\'{e}clet number, $Pe = \text{diam}(\Omega)\|\vb\|/\varepsilon$, is large, the flow is dominated by the convective term, which drives sharp changes in the solution value near the boundaries of $\Omega$.  The presence of these so-called \textit{boundary layers} presents a two-fold challenge in the numerical approximation of $u(x,y)$, namely, in discretizing the PDE in such a way that these layers are correctly represented, no matter how large $Pe$ (or how small $\varepsilon$) is, and in efficiently solving the resulting systems of linear equations.  A key technique in the accurate discretization of such problems is the use of boundary-fitted meshes \cite{stynes2018convection,farrell1995design,linss2009layer} coupled with stable discretizations of the equation to achieve accuracy independent of $\varepsilon$.  In this paper, we focus on the commonly used Shishkin meshes, first in the case of tensor-product grids of the unit square, then in the more complicated case of a non-rectangular domain.  However, the techniques proposed here can easily be extended to more general layer-adapted meshes, such as Bakhvalov-type meshes.

Accurate finite-difference and finite-element discretizations of convection-dominated equations have long been studied.  In the finite-difference case, while centred difference discretizations of convective terms offer second-order accuracy (at least on uniform meshes), they lead to non M-matrix structure of the resulting linear system, which can lead to non-physical oscillations in the numerical solution when the mesh P\'{e}clet number, $Pe_h = h\|\vb\|/\varepsilon$, is large, where $h$ is a representative mesh size (such as the diameter of the largest element in the mesh).  The standard remedy for this is to consider a one-sided upwind discretization of the convective term.  Using the standard first-order upwind discretization restores the M-matrix property of the discretized system and can be used to prove stability (and lack of unphysical oscillations) in the discretized systems~\cite{stynes2018convection}.  Similar difficulties arise in the finite-element case, leading to the use of stabilized finite-element methods \cite{silvester2005finite, linss2009layer}.  In this paper, we use the Streamline Upwind Petrov-Galerkin (SUPG) stabilization, which has been frequently considered for singularly perturbed equations~\cite{LIN20013527, StynesTobiska2003, linss2005anisotropic, FRANZ20081818, augustin2011assessment, MR2869032}.

  The focus of this paper is on the efficient iterative solution of the linear systems that result from the SUPG finite-element discretization of singularly perturbed convection-diffusion equations on Shishkin meshes.  The question of efficient linear solvers for upwind finite-difference discretizations on these meshes is well-studied, starting with the study of diagonal (Jacobi) preconditioning~\cite{hans1996note} and the Gauss-Seidel iteration~\cite{farrell1998convergence} in one dimension, along with incomplete LU preconditioning for the two-dimensional problem~\cite{AnHe03}. There is a growing literature studying Schwarz-type domain decomposition methods for the upwind discretization as well, starting with the study of multiplicative Schwarz with volumetric overlap for the one-dimensional problem in the continuum~\cite{TPMathew_1998a} and discretized~\cite{macmullen_oriordan_shishkin_2002} forms.  Multiplicative Schwarz with minimal overlap has been studied more recently, in both one dimension~\cite{echeverria2018convergence} and in two dimensions~\cite{etna_vol54_pp31-50}.  A more sophisticated Schwarz method, with inexact subdomain solves, was also recently proposed~\cite{SMacLachlan_etal_2021a}.  Finally, a multigrid method for the upwind finite-difference discretization has also been considered~\cite{gaspar2000multigrid, gaspar2002some}.  Here, we extend this multigrid method to the SUPG discretization, on structured multigrid hierarchies on both rectangular and more general domains.  We note that multigrid for SUPG on rectangular domains has been well-studied~\cite{MR1715555, MR2238675, doi:10.1137/060662940, MR2295933, Rees_etal_2011, doi:10.1137/17M1144350}, although primarily on uniform rectangular grids.  The present work is most similar to Ramage's 1999 paper~\cite{MR1715555}, which also considered layer-adapted meshes and robust multigrid methods; however, here, we demonstrate that parameter-robust solvers can be achieved when combining layer-adapted meshes tailored to the singular-perturbation parameter and domain under consideration, even in the case of a non-rectangular domain.

  The remainder of this paper is organized as follows.  In Section~\ref{sec:CDdiscretization}, we present a model problem and discuss its discretization using the SUPG finite-element formulation.  A review of multigrid methods and details on their specialization to these discretizations are presented in Section~\ref{SME}.  Numerical results are presented in Section~\ref{sec:numerical}, including results on tensor-product meshes of the unit square, on layer-adapted meshes for annular domains, and for the well-known ``Hemker Problem''~\cite{hemker1996singularly, hegarty2020numerical}.  Finally, concluding remarks are given in Section~\ref{sec:conclusions}.
  
\section{Convection-Diffusion and its Discretization} \label{sec:CDdiscretization}
We begin by considering the model problem:
\begin{equation} \label{eq_mod}
\mathcal{L}u:=-\varepsilon\Delta u - \vb\cdot \nabla u +cu=f ~~\text{in} ~~ \Omega
\end{equation}
subject to homogeneous boundary conditions 
\begin{align*}
~~ u=0~~\text{on}~~ \Gamma=\partial\Omega 
\end{align*}
where $\varepsilon \in (0,1]$ is a small positive perturbation parameter, with $\vb,~c\geq1$, and $f$ being generic terms that represents sufficiently smooth functions. It is well known that small values of  $\varepsilon$ give rise to exponential layers near any outflow boundary~\cite{roos2008robust,miller1996fitted} (where $\vb\cdot\vn > 0$ for outward unit normal vector $\vn$), and to parabolic layers near characteristic boundaries, where $\vb\cdot\vn = 0$.  In this section, we focus on the case where $\Omega = (0,1)^2$ and the convection coefficient is aligned with the $x$-axis,
\begin{align*}
\vb=(b_1,0)\text{ with }b_1\geq\beta_1 ~~ \text{for~a~positive~constant}~ \beta_1,
\end{align*}
but note that the discretization techniques are equally applicable to more general problems and domains, as considered later on.
To ensure well-posedness, we assume that 
$$c+\frac{1}{2}\nabla\cdot\vb\geq \gamma>0. $$
We note that the case of non-grid-aligned convection, with $\vb = (b_1,b_2)$ for positive $b_1$ and $b_2$ is generally considered to possess a less challenging layer structure, due to the absence of parabolic layers.

The standard weak form for~\eqref{eq_mod} arises by multiplying it by a test function, $v\in H^1_0(\Omega)$, and then integrating over the domain, resulting in the problem of finding $u\in H^1_0(\Omega)$ such that $a(u,v) = \langle f,v \rangle$ for all $v\in H^1_0(\Omega)$, for
$$ a(u,v)= \varepsilon \langle\nabla u,\nabla v\rangle-\langle \vb \cdot \nabla u ,v\rangle + \langle cu, v \rangle, $$
where $\langle u,v\rangle = \int_\Omega uv$ is the standard $L^2(\Omega)$ inner product.  For any finite-dimensional subspace, $V^h \subset H^1_0(\Omega)$, the resulting Galerkin approximation is to find $u^h  \in V^h$ such that $a(u^h,v^h) = \langle f,v^h\rangle$ for all $v^h\in V^h$.

\subsection{Streamline Upwind Petrov Galerkin Discretization}\label{sec:SUPG}
\begin{figure}[t]
\begin{minipage}{0.47\textwidth}
\centering
\includegraphics[trim={16cm 7cm 12cm 7cm},clip,scale=.35]{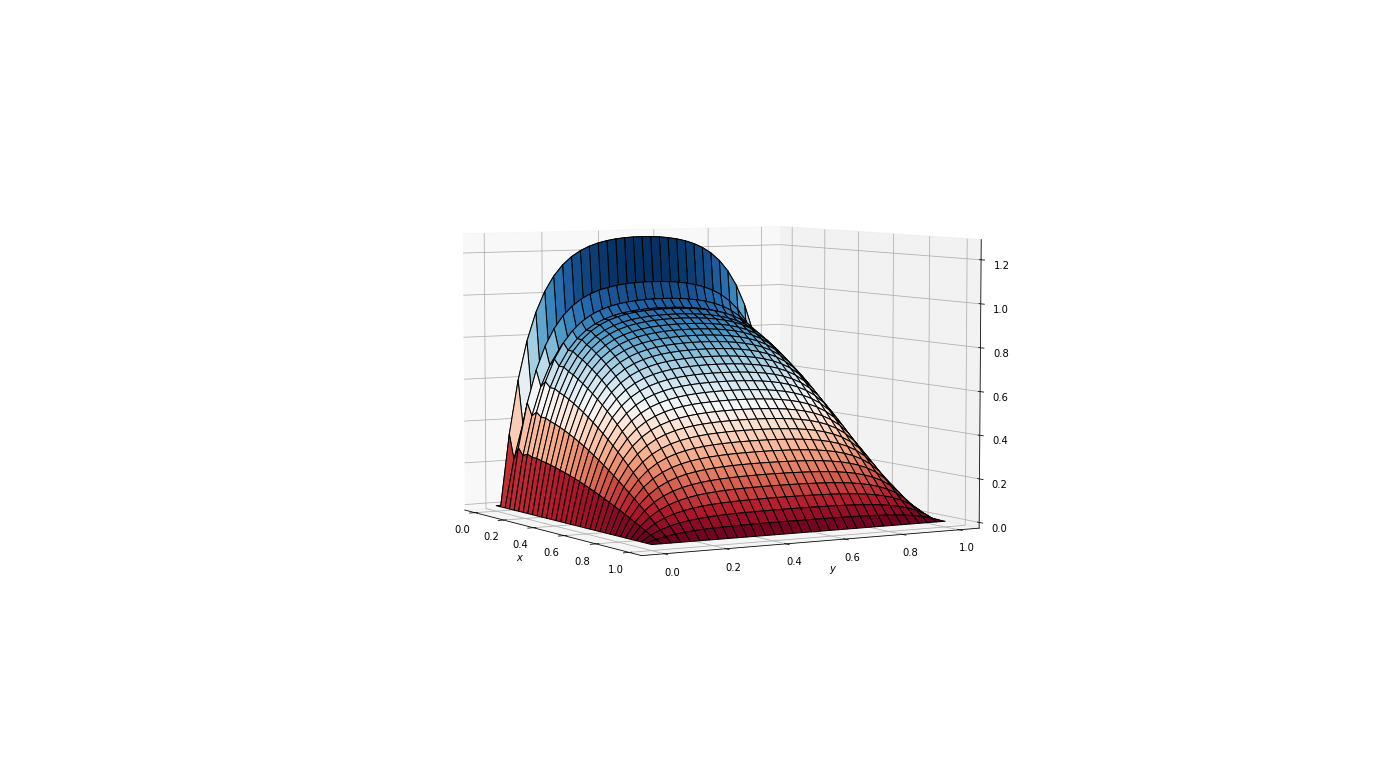}
\end{minipage}
\begin{minipage}{0.47\textwidth}
\centering
\includegraphics[trim={15cm 7cm 12cm 7cm},clip,scale=.35]{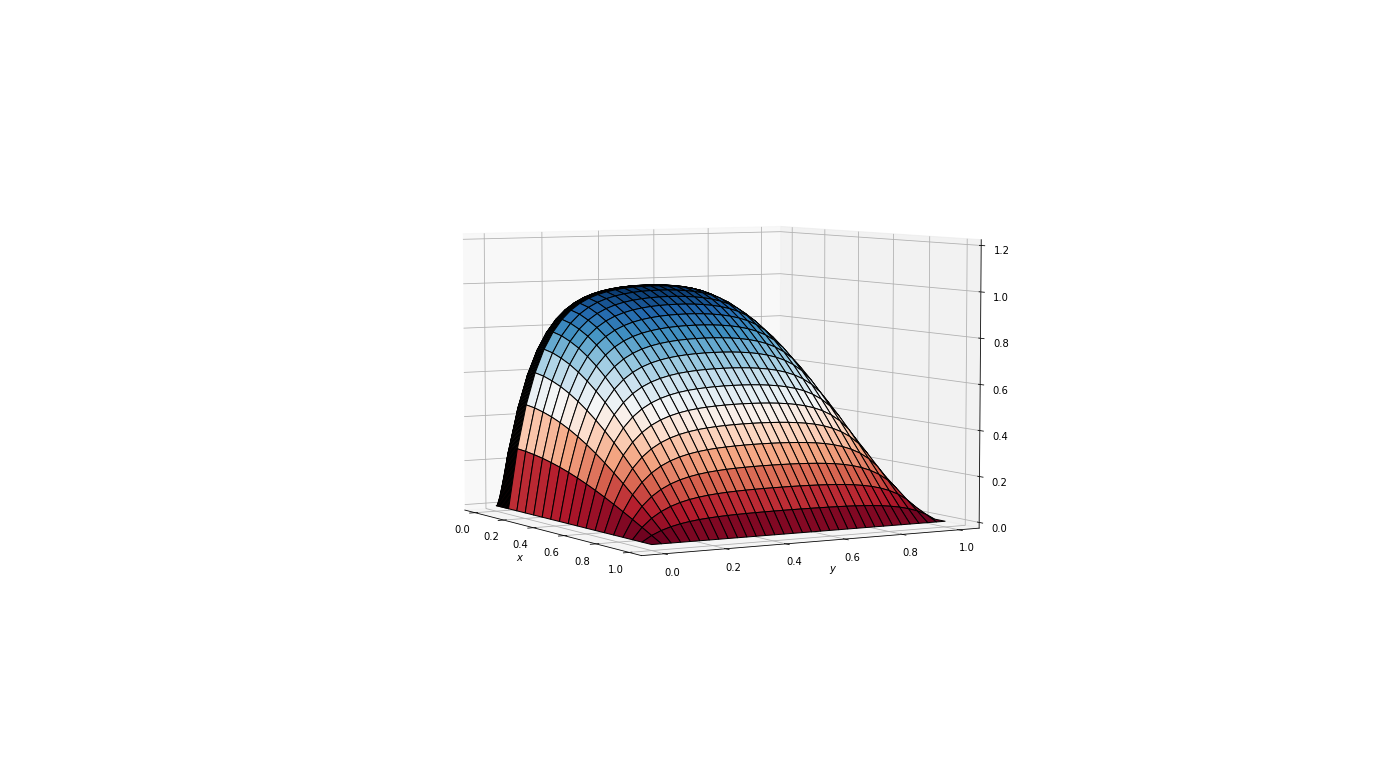}
\end{minipage}
\caption{Model problem solution to Equation~\eqref{eq4} with $\varepsilon=10^{-2}$.  At left, the solution without SUPG shows significant oscillations in the exponential layer.  At right, these wiggles are not present in the solution using SUPG.}
\label{fig21}
\end{figure}

Equation~\eqref{eq_mod} can be categorized on the basis of the mesh P\'{e}clet number, which is the ratio of convective to diffusive transport:
$$Pe_\mathcal{T}=\dfrac{\| \vb \| h_{\mathcal{T}}}{2 \varepsilon}.$$
Here, $\vb$ is the transport velocity and $h_{\mathcal{T}}$  is the diameter of the mesh cell i.e., $h_{\mathcal{T}} = \displaystyle\sup_{\overline{x},x \in \mathcal{T}}\|\overline{x}-x\|$.
If $Pe_\mathcal{T}\leq 1$ for all $\mathcal{T}$ in the mesh, Equation~\eqref{eq_mod} is said to be diffusion dominated and can be accurately solved by the standard Galerkin method. However, when there exist cells for which $Pe_\mathcal{T}\gg 1$, Equation~\eqref{eq_mod} represents a singularly perturbed problem and solution via the standard Galerkin discretization gives rise to spatial oscillations in the approximate solution, commonly known as \emph{wiggles}, as shown at the left of Figure \ref{fig21}. In particular, if a uniform mesh is employed for discretization, then the computed solution is only reliable if an impractically large number of mesh points are used, such as when $h_{\mathcal{T}} \ll \varepsilon$ for all $\mathcal{T}$. To overcome such problems, a common strategy is to use layer-adapted meshes to provide sufficient grid resolution in areas where sharp changes in the solution value are expected~\cite{linss2009layer,miller1996fitted}.  In this paper, we will focus on the so-called ``Shishkin meshes''~\cite{miller1996fitted,LIN20013527}, that are tensor-product meshes with piecewise-constant mesh widths, chosen to reflect the expected layer structure in the solution.  Using this \textit{a priori} knowledge, there are many cases for which we can prove that a numerical approximation, $u^h$, to the true solution, $u$, satisfies
$$\lVert u-u^h\rVert \leqslant C\left(\frac{\log N}{N}\right)^{p}$$
  for positive constants $C$ and $p$ that are independent of $\varepsilon$~\cite{linss2009layer}.

  While the standard Galerkin method on a layer-adapted mesh can yield an adequate approximation to the solution of a singularly perturbed convection-diffusion equation, it is more common to consider more robust discretizations, that directly attempt to avoid numerical oscillations in $u^h$.  Petrov-Galerkin methods for such problems (on uniform grids) were first considered by Christie et al.~\cite{christie1976finite}.
Shortly after this work, Brooks and Hughes proposed the streamline upwind Petrov-Galerkin (SUPG) framework that is now commonplace~\cite{BROOKS1982199}.
 We describe this method beginning with the weak formulation of~\eqref{eq_mod}, to find $u \in V$ such that 
\begin{equation}\label{weak_conv_form}
\varepsilon \langle\nabla u,\nabla v\rangle-\langle\vb \cdot \nabla u ,v\rangle+\langle cu,v\rangle= \langle f,v\rangle \text{ for all }v\in W,
\end{equation}
where we note that we allow $V\neq W$, giving different test and trial spaces.
At the discrete level, we take $V^h$ to be the standard Lagrange finite-element space, but take $W^h$ to be the span of functions $w^h$ given by
\[
w^h =v^h-\tau \vb\cdot\nabla v^h,
\]
for $v^h\in V^h$.  This leads to the weak form of finding $u^h \in V^h$ such that
\begin{equation}\label{eq:SUPG_weak}
a_{S}(u^h,v^h)=f_S(f,v^h)\text{ for all }v^h\in V^h,
\end{equation}
with
\begin{align*}
a_S(u^h,v^h)&:= \varepsilon \langle\nabla u^h,\nabla v^h\rangle-\langle\vb \cdot \nabla u^h ,v^h\rangle+ \langle cu^h ,v^h\rangle- \sum_{\mathcal{T} \in \Omega^h}\tau_\mathcal{T} \langle cu^h-\vb\cdot\nabla u^h-\varepsilon \Delta u^h,\vb\cdot\nabla v^h\rangle_\mathcal{T},\\
f_S(f,v^h)&:=\langle f,v^h\rangle-\sum_{\mathcal{T} \in \Omega^h}\tau_\mathcal{T} \langle f,\vb\cdot\nabla v^h\rangle_\mathcal{T},
\end{align*}
where $\tau_\mathcal{T}$ is the stabilization parameter on element $\mathcal{T}$ for the streamline upwind method, and the summations are over every element in the mesh, $\Omega^h$.
A standard choice of $\tau_\mathcal{T}$ for uniform meshes~\cite{roos2008robust} is given by
\begin{equation}\label{Tau_0}
\tau_{\mathcal{T}}=
    \begin{cases}
      h_{\mathcal{T}}\tau_0 & \text{if } Pe_\mathcal{T}>1\\
      h_{\mathcal{T}}^2\tau_1/\varepsilon &  \text{if } Pe_\mathcal{T}\leq 1\\      
    \end{cases},
\end{equation}
where $Pe_\mathcal{T}$ is mesh P\'{e}clet number, $\tau_0$ and $\tau_1$ are appropriate constants, and $h_{\mathcal{T}}$ is the diameter of the mesh cell.

\begin{figure}[t]
\centering
\includegraphics[width=0.85\textwidth]{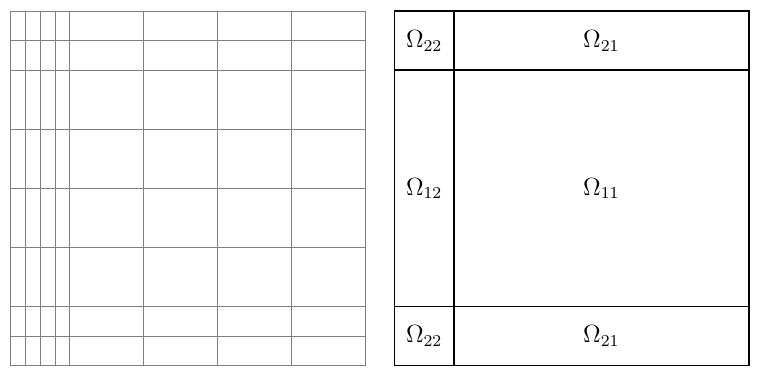}
\caption{Left: Sample Shishkin mesh with $N=8$ for this dissection.  Right: Dissection of $\Omega$ for unit-square domain based on transition points $\lambda_1$ and $\lambda_2$ and expected layer structure.}
\label{fig1}
\end{figure}

We discretize the weak form in~\eqref{eq:SUPG_weak} on quadrilateral Shishkin meshes of the unit square with $N$ mesh intervals in both the $x-$ and $y-$directions, with layered structure chosen to resolve both the \emph{parabolic} boundary layers adjacent to the bottom and top edges ($y=0$ and $y=1$) and the \emph{exponential} boundary layer on the left edge ($x=0$).  To resolve these layers properly, we use a \emph{fitted} mesh that is suitably refined in the boundary regions where these layers occur. This piecewise uniform mesh is given by a tensor-product mesh formed from two piecewise uniform fitted meshes in one dimension. To resolve the exponential layer, we partition the one-dimensional mesh into two pieces based on a transition point
\begin{equation}\label{eq:xShishkin_transition}
\lambda_1:= \min\left\lbrace \frac{1}{2},\sigma\varepsilon \ln N\right\rbrace,
\end{equation}
for suitably chosen parameter, $\sigma$.
Each of the two resulting intervals, $[0,\lambda_1]$ and $[\lambda_1,1]$ are then split into $N/2$ equal mesh elements, for given $N$. 
To resolve the parabolic layers, we define a one-dimensional mesh with two transition points, given by $\lambda_2$ and $1-\lambda_2$, where
\begin{equation}\label{eq:yShishkin_transition}
 \lambda_2:= \min\left\lbrace \frac{1}{4}, \sigma\sqrt{\varepsilon} \ln N\right\rbrace.
 \end{equation}
 The resulting partition of $[0,1]$ into three intervals, $[0,\lambda_2]$, $[\lambda_2,1-\lambda_2]$, and $[\lambda_2,1]$, is then subdivided into $N/4$ mesh elements of equal size for $[0,\lambda_2]$ and $[1-\lambda_2,1]$, while the ``interior'' subdomain $[\lambda_2,1-\lambda_2]$ is partitioned into $N/2$ equal sized mesh elements. Typically, $\sigma$ is taken equal to the order of the method divided by a lower bound on the coefficient $c$, which we will take to be 1.  As we consider piecewise bilinear finite elements on the tensor-product mesh, we take $\sigma = 2.5$.
The transition points divide the unit-square domain, $\Omega$, into 6 subregions, as depicted at right of Figure~\ref{fig1}.  The resulting Shishkin mesh for $N=8$ is depicted at left of Figure~\ref{fig1}.

\section{ Solution Techniques }\label{SME}
Discretization of~\eqref{eq:SUPG_weak} on an $N\times N$ Shishkin mesh results in a linear system to solve that we write as $A_hu_h = b_h$, using a subscript to distinguish between the representation of the finite-element solution as a function, $u^h$, and the vector representing its coefficients in the standard basis, given here as $u_h$.  Direct methods, such as Gaussian elimination, are often the fastest choices for small values of $N$, but are known to scale no better than $\mathcal{O}(N^3)$, using nested dissection ordering~\cite{Ge1973_nested}.  Moreover, such scaling is only known to be achieved on (logically) rectangularly structured meshes in two dimensions, with worse scaling on unstructured meshes in two dimensions, and for all three-dimensional problems.

Thus, for systems such as those resulting from the SUPG discretization of our model problem, iterative methods are more preferable. In contrast to direct methods, iterative methods seek to iteratively improve an approximate solution until an appropriate stopping criterion is satisfied, such as reduction of the norm of the residual below a specified absolute or relative \emph{tolerance}. Iterative methods are broadly classified into \emph{stationary} and \emph{non-stationary} methods. Stationary methods for $A_hu_h=b_h$ are generally written as
\begin{equation*}
 u_h^{m+1} := u_h^m + M_h^{-1} r_h^m,
\end{equation*}
where $r_h^m = b_h - A_hu_h^m$ is the residual after the $m^{th}$ iteration and $M_h^{-1}$ is some computationally feasible approximation of $A_h^{-1}$.  Standard choices of $M_h$ include as the scaled diagonal of $A_h$, leading to the weighted Jacobi iteration, or the lower-triangular part of $A_h$, leading to the Gauss-Seidel iteration, but more sophisticated choices, such as a multigrid V- or W-cycle can also be expressed in this framework.
Non-stationary methods include the well-known class of Krylov subspace methods~\cite{saad2003iterative}.  These methods iteratively construct bases for the Krylov space given by $\text{span}\{r_h^0, A_hr_h^0, A_h^2 r_h^0, \cdots \}$ and look for the \emph{best} solution vector in this space in some measure, with different methods distinguished between the basis constructed and the solution chosen.  For most problems, the most effective solution techniques combine an ``outer'' Krylov subspace algorithm with a preconditioner chosen from a suitable stationary technique.  In this paper, we focus on the use of the flexible variant of GMRES (FGMRES)~\cite{saad1993flexible}, with a multigrid preconditioner.  We note that FGMRES is chosen as the Krylov method for two reasons.  First, it is naturally a right-preconditioned method, meaning that the Krylov method is driven by minimizing the residual norm, $\|b_h - A_hu_h^m\|$, at each iteration, rather than the preconditioned residual norm minimized in left-preconditioned GMRES.  Secondly, FGMRES ``exchanges'' the cost of an extra preconditioner application in classical right-preconditioned GMRES for extra vector storage for the preconditioned Arnoldi vectors.  In our setting, with no serious memory limitations, but a somewhat expensive preconditioner, this is the more efficient choice.

\subsection{Multigrid Methods}
Multigrid methods~\cite{trottenberg2000multigrid, wesseling2004introduction} were originally proposed for the numerical solution of elliptic partial differential equations and have been extended and adapted for a much wider class of problems in recent years.
Here, we briefly introduce multigrid.  Classical stationary iterative methods, such as weighted Jacobi and Gauss-Seidel, are well known to converge slowly for discretizations of simple elliptic problems, such as the Laplacian.  Analyzing such methods generally shows that not all modes of error are reduced at the same rate and, for typical elliptic problems, the error in an approximation can be sensibly partitioned into low- and high-frequency modes, with quick reduction of high-frequency modes and slow reduction of low-frequency modes of the error.  Due to this property, such iterations are known as \emph{smoothing} or \emph{relaxation} methods.

Multigrid gets its name from the fact that the PDE problem is discretized on multiple grids with varying meshsizes~---~usually taken as $h, 2h, 4h, \ldots$ when the finest discretization grid is uniform, although here (in deference to the nonuniform Shishkin meshes) we will index by level instead, with level 1 denoting the finest grid.  Since low-frequency error components are the slowest to be reduced by relaxation on the finest grid, we aim to compute a correction to the approximation produced by relaxation, using a coarser representation of the problem.  Since many low-frequency error components appear as high-frequency components on the coarse grid, simple relaxation there can be more effective at reducing these error components than it is on the finest grid.  This principle is applied recursively to yield a multigrid cycle.  Information is transferred from one grid to the next coarsest in the hierarchy via averaging of the fine-grid residual, referred to as \emph{residual restriction}, and back from a coarse grid to the next finest grid via interpolation of the coarse-grid approximation of the error, known as \emph{error prolongation}.  We note that this prolongation typically both damps low-frequency errors on the fine grid and excites some high-frequency error components, due to frequency aliasing, so relaxation is used again after this \emph{coarse-grid correction} phase to further damp any high-frequency errors on the fine grid.  Thus, we distinguish between \emph{pre-relaxation}, performed before a residual is calculated and restricted to the coarse grid, and \emph{post-relaxation}, performed after the error correction is attained from the coarse grid.  Algorithm~\ref{alg:multigrid} details standard multigrid cycles, with parameters $\nu_1, \nu_2 \geq 0$ to determine the number of pre- and post-relaxation iterations used, and $\gamma_1$ and $\gamma_2$ used to determine the cycle type.  With $\gamma_1 = \gamma_2 = 1$, only a single recursion is performed on each level leading to a V-cycle, while $\gamma_1 = \gamma_2 = 2$ performs a ``double'' recursion on each level, leading to a W-cycle.  The common case of a multigrid F-cycle corresponds to $\gamma_1 = 1$ and $\gamma_2 = 2$, where we recurse one time on the ``downward'' sweep (fine-to-coarse) of the multigrid cycle, and again on the ``upward'' sweep.

As is typical in multigrid, Algorithm~\ref{alg:multigrid} is composed of two pieces, a ``setup'' phase, where the grid hierarchy itself and operators on all grids in the hierarchy are pre-computed (possibly including information needed for the relaxation sweeps on these grids as well), and a ``solve'' phase, that includes pre- and post-relaxation on each grid, along with the recursive calls to the next coarser grid in the hierarchy.  We write the algorithm from the perspective of level $l$ in a hierarchy of $C$ levels, with $1 \leq l \leq C$, where the number of levels is known beforehand.  For the solve phase, the recursion is clear; the setup phase is called independently on each level, with the assumption that $A_l$ is known (since $A_1$, the discretization matrix on the finest level, is known at the start of the coarsening process).
 The (residual) restriction operator is denoted by $I_l^{l+1}$, and the (error) prolongation operator by $I_{l+1}^l$.

\begin{algorithm} 
\centering {$u_l^{m+1} = {\rm \textbf{MG}}(l, u_{l}^m,b_l,\nu_1,\nu_2, \gamma_1, \gamma_2)$.}
\caption{\bf ~~ Multigrid pseudocode}
\begin{enumerate} \sf 
\item[{\bf (0)}] \sf  {\bf Setup Phase}\\
  -- \sf \underline{If} $l=C$, compute $A_l^{-1}$\\
  -- \sf \underline{Else} Establish static access to:\vspace{-2mm}
  \begin{itemize}
    \item ${\rm \textbf{relax}}(A_l, u^m_l, b_l)$ (relaxation on level $l$).
	 \item $I_l^{l+1}$ (the restriction operator from level $l$ to level $l+1$).
	 \item $I_{l+1}^l$ (the prolongation operator from level $l+1$ to level $l$).
	 \item $A_{l+1}$ (the coarse-grid operator on level $l+1$).
	\end{itemize}
\item[{\bf (1)}] {\bf Solve Phase}
\begin{itemize}
\item[{\bf (a)}] {\bf Pre-relaxation} \\[1.0ex]
  -- \sf  Compute
$\overline{u}^m_l$  by applying $\nu_1$ relaxation
steps to $u^m_l$:
$  
\overline{u}^m_l = {\rm \textbf{relax}}^{\normalsize\mbox{$\nu_1$}} ( A_l, u^m_l,
b_l)\enspace.
$
\item[{\bf (b)}] \sf  {\bf Coarse grid correction} \\[1.0ex]
\begin{tabular}{lll}
-- &   Compute the residual  & $  \overline{r}^m_l = b_l - A_l \overline{u}^m_l$ \enspace.\\[1.0ex]
-- &   Restrict the residual & $  \overline{r}^m_{l+1} = I_l^{l+1} \hspace{2mm} \overline{r}^m_l$ \enspace.\\[1.0ex]
-- &   Compute the approximate error   $\widehat{e}^m_{l+1}$ from the \textit{defect equation}. & $A_{l+1} \hspace{2mm} \widehat{e}^m_{l+1} = \overline{r}^m_{l+1}$\\
   &	by the following mini algorithm& \\
\end{tabular}

\hspace*{5ex}\framebox[10cm]{
\begin{minipage}{10cm}
\vspace{2mm}
 \sf 
\underline{If} $l+1=C$, $\hat{e}_{l+1}^{m} = A_{l+1}^{-1}\overline{r}_{l+1}^m$\\
\underline{Else} approximate $\widehat{e}_{l+1}^{m}$ iteratively by the recursion:\\
$\text{$\hspace{3mm}$}$ $\widehat{e}_{l+1}^{m,1} = \overline{0}$; \\
$\text{$\hspace{3mm}$}$\underline{do} $i = 1$ to $\gamma_2$  \\
$\text{$\hspace{5mm}$}$ \underline{If} $\hspace{1mm} i = 1$,\\
$ \text{$\hspace{10mm}$} \widehat{e}_{l+1}^{m,i+1}\hspace{-1mm}={\rm\textbf{MG}}(l+1,\widehat{e}_{l+1}^{m,i},\overline{r}_{l+1}^m,\nu_1,\nu_2,\gamma_1,\gamma_2)$\\
$\text{$\hspace{5mm}$}$ \underline{Else} \\
$\text{$\hspace{10mm}$} \widehat{e}_{l+1}^{m,i+1}\hspace{-1mm}={\rm\textbf{MG}}(l+1,\widehat{e}_{l+1}^{m,i},\overline{r}_{l+1}^m,\nu_1,\nu_2,\gamma_2,\gamma_1)$\\
\end{minipage}}
\vspace*{2mm}

\begin{tabular}{lll}
-- & \sf  Interpolate the correction \hspace*{8ex} 
      &$  \widehat{e}^m_l = I_{l+1}^l \hspace{2mm} \widehat{e}^m_{l+1}$ \enspace.\\[1.0ex]
-- & \sf  Compute the corrected  approximation on $\Omega_l$ 
                & $  \widehat{u}_l^{m} = \overline{u}^m_l + \widehat{e}^m_l$\enspace.
\end{tabular}
\\
\item[{\bf (c)}]  \sf {\bf Post--relaxation}
-- \sf  Compute $u^{m+1}_l$ by applying $\nu_2$ relaxation
        steps to $\widehat{u}_l^{m}$:
$   u^{m+1}_l = {\rm \textbf{relax}} ^{\normalsize\mbox{$\nu_2$}} 
    (A_l, \widehat{u}_l^{m} , b_l) \enspace.
        $
\end{itemize}
\end{enumerate} 
 \label{alg:multigrid}
\end{algorithm}

For standard elliptic problems discretized on equidistant meshes, multigrid is well-known to work efficiently. However, preliminary experiments showed that its \emph{components} need to be adjusted for the case of interest here, where mesh sizes vary in each spatial dimension. A closely related case is that of an anisotropic diffusion equation, where it is well-known that either adjustment in the \emph{grid coarsening scheme} (to coarsen only in the direction of strong connections) or in the \emph{relaxation scheme} or both is required to maintain robust convergence. To handle anisotropy in the relaxation scheme, block-relaxation methods that simultaneously relax strongly connected unknowns are commonly used~\cite{brandt1977multi,trottenberg2000multigrid}.  As expected (and observed in our own experiments), the convergence of standard multigrid methods suffer in cases where mesh aspect ratios are much over or under one, as is the case with the current model problems on Shishkin meshes~\cite{gaspar2002some,gaspar2000multigrid}.  Thus, we next consider how to construct robust block-relaxation methods for the problem at hand.

\subsection{Relaxation Scheme}\label{subsec:RELSCH}
A successful multigrid scheme fundamentally depends on the complementarity (in error reduction) between its two sub-procedures, \emph{relaxation}, and \emph{coarse-grid correction}. In the current scenario, with anisotropy induced mainly due to discretization of the PDE on Shishkin meshes, this complementarity may be ensured by the use of \emph{block relaxation} and \emph{standard coarsening} as our multigrid method's main components. 
The particular block relaxation that we use is based on grouping spatial lines of grid-points as \emph{blocks} to be relaxed simultaneously, i.e., \emph{Line Relaxation} that corresponds to vertical / horizontal lines.  Typical descriptions of line relaxation focus on variants of block Gauss-Seidel (see, for example, Trottenberg et al.~\cite{trottenberg2000multigrid}), including its parallelizable ``red-black'' (often called ``zebra'') variant, and variants that alternate between lines parallel to the $x$-axis and those parallel to the $y$-axis.  Here, we follow an additive variant, where a single residual is computed, and updates are computed for a specified set of lines from that common residual.  These updates are summed to form the correction to the approximation that generated the initial residual.  In order to ensure proper relaxation (in the multigrid sense), we must appropriately weight this correction, which we do either by using a Chebyshev iteration on each level, with the interval chosen (roughly) to damp the upper 3/4 of the spectrum of the line-relaxation-preconditioned linear system (since we are using standard coarsening in two dimensions), or simply by using GMRES to automatically weight the correction.  Details of these are discussed in the numerical results below.

While specialized routines to implement block relaxation on data structures tailored to fully structured grids are somewhat straightforward to implement, robust and general block relaxation schemes that can be used with general-purpose discretization codes are not.  Here, we make use of and extend the PCPatch~\cite{farrell2021pcpatch} implementation in Firedrake~\cite{FiredrakeUserManual, kirby2018solver} and PETSc~\cite{balay2019petsc} to allow general construction of block relaxation schemes based on topological properties of the underlying mesh data structure (here, in particular, the node locations).  A critical part of this is the specification of which discrete DoFs belong to each block (or patch, in the PCPatch terminology).  The original PCPatch implementation provided a basic mechanism for line or plane relaxation, allowing the user to specify one or more axes to indicate the relaxation direction, along with a number of divisions and direction for each axis.  For each specified axis, the points in the discretization mesh are sorted by their coordinate value along that axis (so, for the $x$-axis, they are sorted by their $x$-coordinates, leading to lines of points in the $y$-direction for a two-dimensional mesh), and divided into the specified number of evenly-spaced divisions (by the coordinate along the specified axis).  For additive sweeps, the direction plays no important role, but it specifies the direction of sweep for the multiplicative case (allowing ``forward'' and ``backward'' sweeps along each axis).  The PCPatch abstraction then gathers PETSc ``index sets'' of the DoFs for each block specified in this way, and those index sets are used to define the block matrices for the relaxation.

While this allows for proper treatment of patches in the case of anisotropic diffusion operators discretized on uniform rectangular meshes, it does not allow for the block relaxation construction needed on Shishkin meshes, even for rectangular domains.  Thus, we have implemented two extensions to the line/plane relaxation implementation in Firedrake.  First, rather than assuming evenly-spaced divisions, we now allow the user to pass an array of divisions, to allow sensible construction of line relaxation on non-uniform tensor-product meshes.  Secondly, instead of restricting the construction of blocks to be based upon a single coordinate axis, we now allow the user to provide a \textit{key} function that maps the Cartesian coordinates of each grid point to a single real value.  These key values are then sorted and divided into block relaxation sets based on either a uniform division, or a provided (nonuniform) array of divisions.  This functionality is critical in allowing us to define sensible block relaxation schemes on non-Cartesian meshes, such as for the Hemker problem~\cite{hemker1996singularly, hegarty2020numerical}. Both extensions are now available in Firedrake.

To complete this section, we provide an algorithmic sketch of the setup and iteration phase of the block relaxation procedure in Algorithm~\ref{alg:blkGS}.  For the setup phase of this algorithm, we use the system matrix, $A_l$, on level $l$ of the hierarchy, along with a matrix $B_l$, used to describe the block structure.  Matrix $B_l$ has the same number of rows as $A_l$, but a number of columns, $m_l$, equal to the number of blocks to be used on level $l$.  Each column of $B_l$ represents a single block with nonzero entries in rows of $B_l$ indicating membership in that relaxation block.  While we write the algorithm assuming $B_l$ is provided as a binary matrix, we note that this is equivalent to simply providing a list (or set) of integer DoF indices for each block in the relaxation.  This is, in practice, the data structure used in PCPatch, where each block is represented by a PETSc Index Set.  For each block, the block sub-matrix is extracted based on these indices, and the LU factorization of each block sub-matrix is pre-computed.  During the iteration phase, we loop over the blocks, additively accumulating the correction from each block.  While we write this for an arbitrary number of sweeps of relaxation, $\nu$, we note that our typical relaxation scheme is to call this for a single iteration at a time, and to combine multiple iterations using a Chebyshev or GMRES outer iteration.  In this mode, each time we call the relaxation routine, $b_l$ is the current Chebyshev residual or GMRES Arnoldi vector, and $u_l^m$ is a zero vector, and the return vector, $u_l^{m+1}$, is the so-called ``preconditioned residual'' used in those algorithms.  We note that the only change for a multiplicative scheme is to update the global residual inside the loop over $j$, rather than outside this loop.

\begin{algorithm}[!htb] 
\centering {$u_l^{m+1} = {\rm \textbf{BLKRELAX}}(A_l, B_l, u_{l}^m, b_l, \nu)$.}
\caption{\bf ~~ Block Relaxation pseudocode}
\begin{enumerate} \sf 
\item[{\bf (0)}] \sf  {\bf Setup Phase} \\[1.0ex]
  --  \sf Establish static access to the matrix $A_l$ of order $n_l \times n_l$.\\
  --  \sf Establish static access to the block matrix, $B_l$, of order $n_l \times m_l$.\\
  --  \sf \underline{do} $j = 1$ to $m_l$\\
  \begin{itemize}
  \item \sf Extract submatrix $A_l^{(j)}$ corresponding to indices in column $j$ of $B_l$.
    \item \sf Compute LU factorization of $A_l^{(j)}$.
  \end{itemize}
\item[{\bf (1)}] {\bf Iteration phase.} 
  \begin{itemize}
    \item[--] \sf Set $u = u_l^m$
\item[--] \sf  \underline{do} $i = 1$ to $\nu$  \\
  \begin{itemize}
    \item[--] \sf $r = b_l - A_l u$
    \item[--] \sf  \underline{do} $j = 1$ to $m_l$  \\
      \begin{tabular}{lll}
         -- & Restrict residual to block $j$ & $r^{(j)} = r(B_l(:,j))$\\
         -- & Solve block $j$ subsystem & $A_l^{(j)}u^{(j)} = r^{(j)}$\\
         -- & Update approximation & $u(B_l(:,j)) = u(B_l(:,j)) + u^{(j)}$\\
      \end{tabular}
    \item[--] \sf \underline{continue} j  
  \end{itemize}          
\item[--] \sf \underline{continue} i
  \item[--] \sf $u_l^{m+1} = u$.
  \end{itemize}
  \end{enumerate}
 \label{alg:blkGS}
\end{algorithm}

\begin{figure}
\begin{minipage}{.47\textwidth}
\centering
\includegraphics[scale=.3]{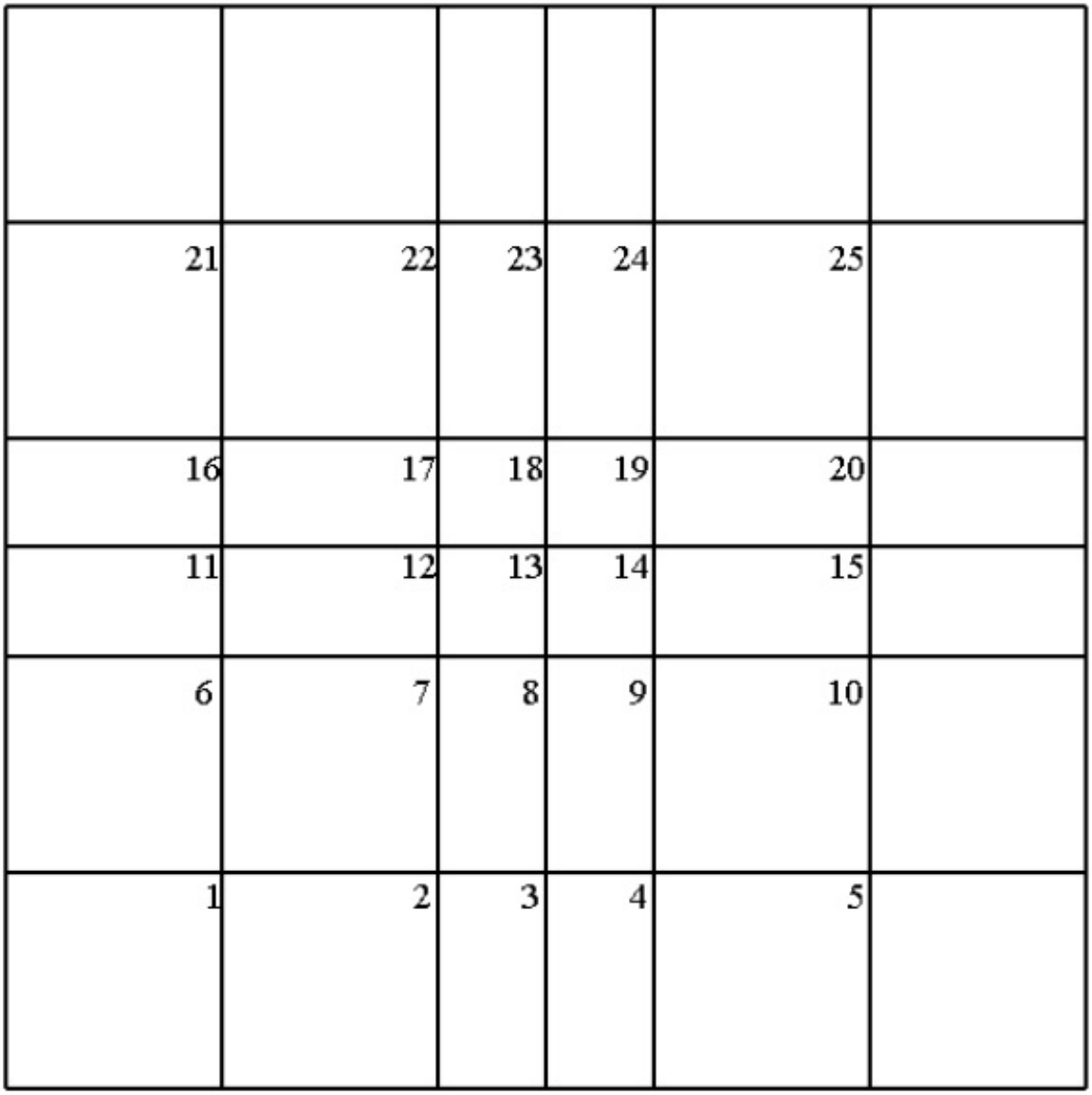}
\caption{Structured enumeration}
\end{minipage}
\begin{minipage}{.47\textwidth}
\centering
\includegraphics[scale=.3]{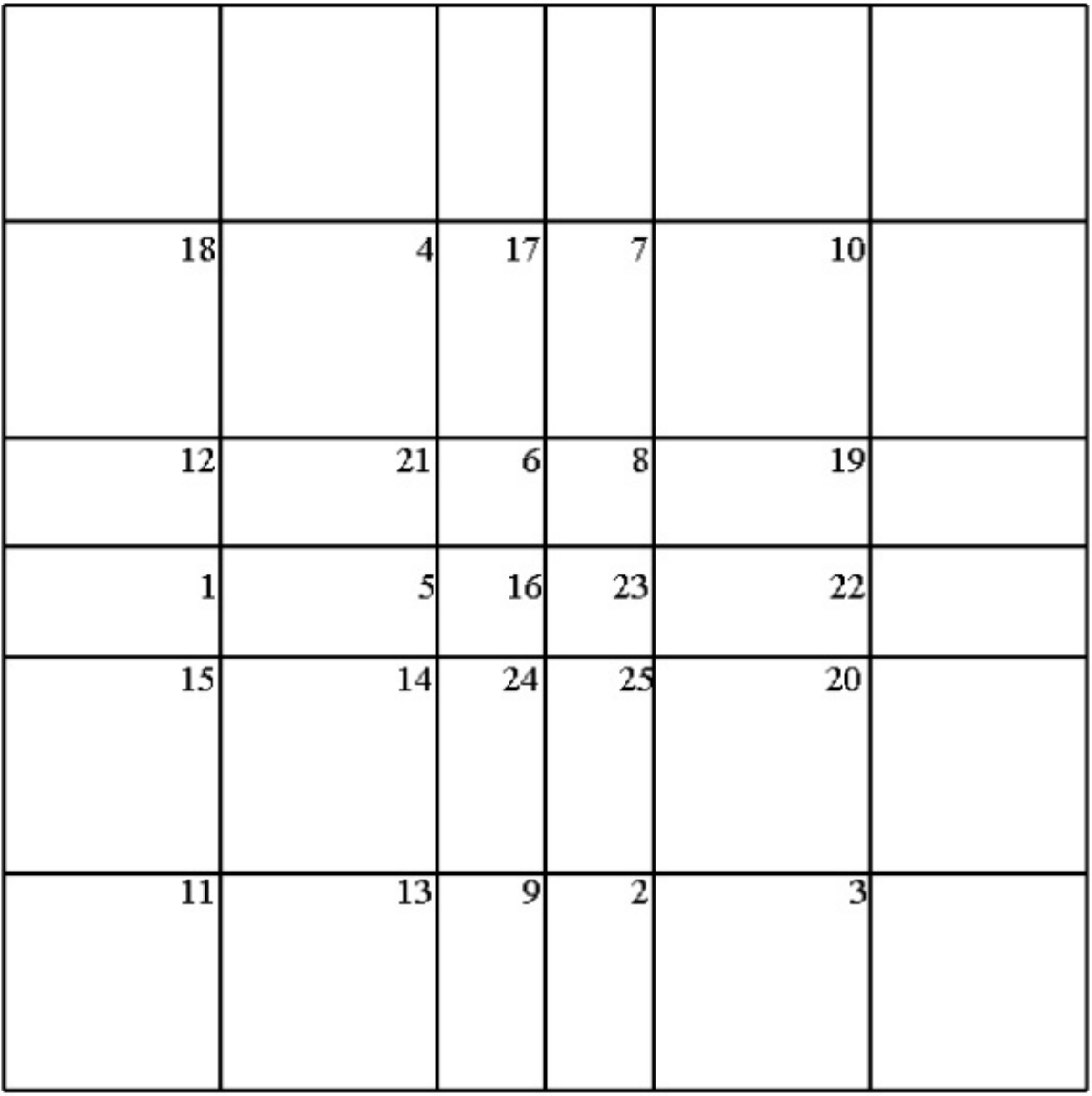}
\caption{Unstructured enumeration}
\end{minipage}
\caption{Examples of rectangular tensor-product grids with structured and unstructured enumeration.}
\label{fig:2dgrids}
\end{figure}

To illustrate the structure expected in $B_l$, we consider two examples of possible enumeration of a structured grid of a two-dimensional rectangle.
\begin{example}\label{exm1:structured}
  Consider the rectangular grid shown at left of Figure \ref{fig:2dgrids}, with a \emph{lexicographic} enumeration.  The ordering of the enumeration naturally suggests lines in the $x$-direction (orthogonal to the $y$-axis in the PCPatch specification), that are captured by the matrix
  \[
B_l^T = \left [ \begin{array}{ccccccccccccccccccccccccc}			    
                            1&1&1&1&1&0&0&0&0&0&0&0&0&0&0&0&0&0&0&0&0&0&0&0&0 \\
                            0&0&0&0&0&1&1&1&1&1&0&0&0&0&0&0&0&0&0&0&0&0&0&0&0 \\
                            0&0&0&0&0&0&0&0&0&0&1&1&1&1&1&0&0&0&0&0&0&0&0&0&0 \\
                            0&0&0&0&0&0&0&0&0&0&0&0&0&0&0&1&1&1&1&1&0&0&0&0&0 \\
                            0&0&0&0&0&0&0&0&0&0&0&0&0&0&0&0&0&0&0&0&1&1&1&1&1 \\                            
                            \end{array}
                     \right ],
\]
where we present the transpose of $B_l$ simply to save space.
Lines in the opposite direction would have a corresponding change in the nonzero pattern in $B_l$, with
\[
B_l^T = \left [ \begin{array}{cccccccccccccccccccccccccc}
                            1&0&0&0&0&1&0&0&0&0&1&0&0&0&0&1&0&0&0&0&1&0&0&0&0 \\
                            0&1&0&0&0&0&1&0&0&0&0&1&0&0&0&0&1&0&0&0&0&1&0&0&0 \\
                            0&0&1&0&0&0&0&1&0&0&0&0&1&0&0&0&0&1&0&0&0&0&1&0&0 \\
                            0&0&0&1&0&0&0&0&1&0&0&0&0&1&0&0&0&0&1&0&0&0&0&1&0 \\
                            0&0&0&0&1&0&0&0&0&1&0&0&0&0&1&0&0&0&0&1&0&0&0&0&1 \\
                            \end{array}
                     \right ].
\]
Line relaxation in both directions would be represented by simply concatenating these two matrices.  Alternative relaxation schemes, such as red-black (zebra) relaxation along lines parallel to the $y$-axis can also be expressed in this notation, as
\[
B_l^T = \left [\begin{array}{cccccccccccccccccccccccccc}
                            1&0&1&0&1&1&0&1&0&1&1&0&1&0&1&1&0&1&0&1&1&0&1&0&1 \\
                            0&1&0&1&0&0&1&0&1&0&0&1&0&1&0&0&1&0&1&0&0&1&0&1&0 \\
                            \end{array}
                     \right ].
\]
\end{example}

\begin{example}\label{exm2:unstructured}
  An alternative numbering for the rectangular grid is shown at right of Figure~\ref{fig:2dgrids}.  We include this example to emphasize that the construction of block relaxation schemes is not limited to lexicographic grid enumeration, and is compatible with any numbering scheme (including the random choice at right of Figure~\ref{fig:2dgrids}).  In this numbering, the columns of $B_l$ are no longer so nicely structured, but the same relaxation schemes can be expressed.  For example, relaxation along lines parallel to the $y$-axis is specified by
  \[
B_l^T = \left [ \begin{array}{ccccccccccccccccccccccccc}			    
                            1&0&0&0&0&0&0&0&0&0&1&1&0&0&1&0&0&1&0&0&0&0&0&0&0 \\
                            0&0&0&1&1&0&0&0&0&0&0&0&1&1&0&0&0&0&0&0&1&0&0&0&0 \\
                            0&0&0&0&0&1&0&0&1&0&0&0&0&0&0&1&1&0&0&0&0&0&0&1&0 \\
                            0&1&0&0&0&0&1&1&0&0&0&0&0&0&0&0&0&0&0&0&0&0&1&0&1 \\
                            0&0&1&0&0&0&0&0&0&1&0&0&0&0&0&0&0&0&1&1&0&1&0&0&0 \\                                                        \end{array}
                     \right ],
\]
while zebra line relaxation parallel to the $x$-axis becomes
\[
B_l^T = \left [ \begin{array}{cccccccccccccccccccccccccc}
                            1&1&1&1&1&0&1&0&1&1&1&0&1&0&0&1&1&1&0&0&0&1&1&0&0 \\
                            0&0&0&0&0&1&0&1&0&0&0&1&0&1&1&0&0&0&1&1&1&0&0&1&1 \\
                            \end{array}
                     \right ].
\]
\end{example}

\subsection{Identifying relaxation lines}\label{ssec:identifying}
The classical wisdom of line relaxation in multigrid methods is that line relaxation should be done in the direction of strong connections in the mesh (see, for example, \cite[Section 5.1.3]{trottenberg2000multigrid}.  Considering the discretization of~\eqref{eq:SUPG_weak}, we can identify the strong and weak connections at a node as a function of the values of the mesh sizes adjacent to the node and the parameter $\tau_{\mathcal{T}}$.  We start by considering the case of tensor-product meshes of rectangular domains.

Using bilinear basis functions, we can express the matrices associated with the individual terms in~\eqref{eq:SUPG_weak} on rectangular meshes as Kronecker products of the one-dimensional discretization stencils.  Given nodes $x_{i-1}$, $x_i$, and $x_{i+1}$, with $h^{(1)} = x_i-x_{i-1}$ and $h^{(2)} = x_{i+1}-x_i$, we can directly compute these stencils as
\begin{align*}
  \langle u',v' \rangle & \rightarrow \left[ \frac{-1}{h^{(1)}} , \frac{1}{h^{(1)}} + \frac{1}{h^{(2)}} , \frac{-1}{h^{(2)}} \right], \\
  \langle u', v \rangle & \rightarrow \left[ \frac{1}{2}, 0, \frac{-1}{2} \right], \\
  \langle u, v \rangle & \rightarrow \left[ \frac{h^{(1)}}{6}, \frac{h^{(1)}}{3} + \frac{h^{(2)}}{3}, \frac{h^{(1)}}{6} \right].
\end{align*}
If we write the matrices associated with these three stencils as $L$ (for Laplacian), $D$ (for first-derivative), and $M$ (for the mass matrix), then, in the constant coefficient case, the first three terms in~\eqref{eq:SUPG_weak} lead to a stiffness matrix of the form
\[
\varepsilon\left(L_x \otimes M_y + M_x \otimes L_y\right) - b_1 D_x \otimes M_y + c M_x \otimes M_y,
\]
where we use the subscript $x$ or $y$ to denote a matrix formed over the mesh in the $x$ or $y$ direction, respectively.
A key realization in determining strong connections in the matrix is the realization that, when $h$ is the maximum mesh size in both directions, the entries in $K_x \otimes M_y$ are $\mathcal{O}(h)$ and those in $M_x \otimes M_y$ are $\mathcal{O}(h^2)$.  Thus, the only possibility for $\mathcal{O}(1)$ entries in these terms come from the diffusive term, which contains factors of $h_x/h_y$ and $h_y/h_x$ (from the two Kronecker product terms) multiplied by $\varepsilon$.  On a typical Shishkin mesh for these problems, we have mesh widths of $\mathcal{O}(\varepsilon/N)$ and $\mathcal{O}(1/N)$ in the $x$-direction and $\mathcal{O}(\sqrt{\varepsilon}/N)$ and $\mathcal{O}(1/N)$ in the $y$-direction, so the ratio $h_x/h_y$ can be $\varepsilon$ (when $x$ is refined, but $y$ is not), $\sqrt{\varepsilon}$ (when $x$ and $y$ are both refined), $1$ (when neither is refined), or $1/\sqrt{\varepsilon}$ (when $y$ is refined, but $x$ is not).  Clearly the ratio $h_y/h_x$ takes the inverse values, ranging from $\sqrt{\varepsilon}$ to $1/\varepsilon$.  Thus, from these terms, the strongest coupling is in the $x$-direction, with $\mathcal{O}(1)$ entries in the matrix in the region $\Omega_{12}$ (see Figure~\ref{fig1}).  The coupling in the $y$-direction can also be strong, with $\mathcal{O}(\sqrt{\varepsilon})$ entries in $\Omega_{21}$, compared to $\mathcal{O}(1/N)$ entries in the $x$-direction in this region.  We will explore if this coupling is strong enough to benefit from line relaxation along lines parallel to the $y$-axis in the numerical results below.

The analysis above disregards two important points.  First, the contributions from the SUPG stabilization terms have not been accounted for.  This is because we are careful to use meshes that resolve boundary layers present in the solution, so we are able to take $\tau_{\mathcal{T}} = \mathcal{O}(h_{\mathcal{T}})$ or smaller.  Consequently, the SUPG terms do not contribute to the strong connections in the matrix.  Secondly, the analysis above only considers the case of tensor-product meshes of rectangular domains.  However, we note that many of the conclusions about scaling of matrix entries remains true, even on the curvilinear meshes considered below, since the mass matrix always has entries that scale with the area of two-dimensional elements, and the convection term will also scale like $h$, leaving the diffusion term as the only possible contributor to strong connections in the matrix, and those strong connections will be aligned with directions in which we are using refined meshes.  We note that the literature on solvers for singularly perturbed convection-diffusion equations primarily focuses on the case where both components of $\vb$ are nonzero and the mesh adaptation is uniform in both directions~\cite{ramage1999multigrid,gaspar2002some}, and the arguments above justify the use of line relaxation in both directions (particularly when Shishkin meshes are used to resolve the two exponential layers typical of outflow boundaries in this case).  The numerical results below will demonstrate that this is not needed for the case of grid-aligned convection with meshes chosen to resolve the resulting parabolic layers, since the $\mathcal{O}(\sqrt{\varepsilon})$ connections in the matrix induced by such layers remain weak precisely in the small-$\varepsilon$ case of interest.

\section{Numerical Results} \label{sec:numerical}
Numerical experiments are performed using the Python interface to Firedrake~\cite{FiredrakeUserManual}, an open-source finite-element package that offers a close interface with PETSc~\cite{balay2019petsc, kirby2018solver}.  As described above, we use PCPatch~\cite{farrell2021pcpatch} to realize the relaxation scheme proposed herein, with the corresponding extensions in its Firedrake interface discussed above.  The codes and major components of Firedrake needed to reproduce the numerical experiments are archived on Zenodo~\cite{zenodo/Firedrake-20231206.0}.
For the first two model problems, we use rectangular meshes with quadrilateral elements, specifying node locations to match those of the Shishkin meshes described above.  For these problems, we generate the coarsest grid in the multigrid hierarchy as an $8\times 8$ mesh, with transition points chosen to match those desired on the finest grid, and subsequent meshes are generated using uniform refinement to the desired resolution.  This procedure is generalized with mesh mapping for the semi-structured grids of the annulus considered next.
For the Hemker problem, described last, we generate the coarsest grid using gmsh~\cite{geuzaine2009gmsh}, create a hierarchy of uniformly refined meshes of this coarse representation, then map these meshes to those obtained by uniform refinement in $(r,\theta)$ on the curved part of the domain.  In all cases, we use a bilinear, $Q_1(\Omega^h)$, approximation space, unless otherwise noted.  As the choice of stabilization parameter, $\tau_{\mathcal{T}}$ in~\eqref{eq:SUPG_weak} is problem-dependent, we describe the choices used with each problem below.

For all problems, we consider solution using FGMRES as the outer Krylov iteration, with the block-relaxation-multigrid preconditioner.  In all cases, we use multigrid V-cycles as the main iteration, with variation in the number of relaxation steps as described below.  On the coarsest grid of the hierarchy, we use a direct solver, via the PETSc interface to MUMPS~\cite{amestoy2000mumps}.  While it is possible to reuse some of the solver infrastructure between solves, particularly on regular grids of the same dimension, we time each experiment for a ``complete'' setup and solve phase, with no such reuse.  Computational times (in minutes) are reported for serial execution on a machine with two 8-core 1.7GHz Intel Xeon Bronze 3106 CPUs and 384 GB of RAM.
For the computed solutions, we measure errors in 3 ways relative to a given analytical or reference solution.  First, we consider the pointwise error between the analytical solution evaluated on the gridpoints of the mesh and the nodal values of the finite-element solution, denoted by ``Max Errors'' in the tables.  
Denoting the standard $H^1(\Omega)$-seminorm and $L^2(\Omega)$ norms by $\mid\cdot\mid_{1,\Omega}$ and $\parallel\cdot\parallel_{\Omega}$, respectively, we define the energy norm by
\begin{equation*}
\Vert v \Vert ^2_\varepsilon := \varepsilon|v|^2_{1,\Omega}+\Vert v \Vert_\Omega^2,
\end{equation*}
and report errors in this norm as ``Energy Errors'' in the tables.  Finally, we introduce the streamline diffusion norm associated with $a_S(\cdot,\cdot)$, defined on the trial function space as
\begin{equation}
\Vert v\Vert_{\text{SD}}=\left( \Vert v \Vert^2_\varepsilon +\sum_{\mathcal{T}\in\Omega^h}\tau_{\mathcal{T}}\Vert \vb \cdot \nabla v \Vert^2_{\mathcal{T}} \right)^{1/2}, \label{eq:SDError}
\end{equation}
and report corresponding ``SD Errors'' in the tables.  We note that the SD norm is often considered a more meaningful measure of solution quality, due to its strong coercivity bound for the SUPG formulation~\cite{StynesTobiska2003, linss2009layer}.

\subsection{Unit Square Test Problems}
The first two test problems considered are taken from Franz et al.~\cite{FRANZ20081818}.

First, we consider the equation
\begin{equation}
-\varepsilon\Delta u - (2-x)u_x+\frac{3}{2}u=f ~~\text{in} ~~ \Omega=(0,1)^2,~~ u=0~~\text{on}~~ \partial\Omega\label{eq4}
\end{equation} 
with the prescribed solution function:
\begin{equation*}
 u(x,y)= \left( \cos \frac{\pi x}{2}-\dfrac{e^{-x/\varepsilon}-e^{-1/\varepsilon}}{1-e^{-1/\varepsilon}}\right)\dfrac{(1-e^{-y/\sqrt{\varepsilon}})(1-e^{-(1-y)/\sqrt{\varepsilon}})}{1-e^{-1/\sqrt{\varepsilon}}}
\end{equation*}
and, consequently, a specific right-hand side, $f$.  We choose the stabilization parameter in~\eqref{eq:SUPG_weak} to be 0 in the left boundary layer, where the mesh size in $x$ is $\mathcal{O}((\varepsilon\ln N)/N)$ and $h_{\mathcal{T}}/2$ in the remainder of the mesh, corresponding to $\tau_0 = 1/2$ and $\tau_1 = 0$ in~\eqref{Tau_0}.

Table~\ref{tab:problem1} presents numerical results for the problem in~\eqref{eq4}, using an outer FGMRES stopping criterion of reducing the absolute residual norm below $1/N^2$ for a problem on an $N\times N$ mesh.  Here, we use V(2,2) cycles, with a Chebyshev iteration for relaxation on each level.  While the interval defining the Chebyshev polynomials could be adjusted by hand, we use a reasonable rule-of-thumb instead, using PETSc's default estimate of the spectrum of the preconditioned operator (computed using preconditioned Arnoldi and a noisy right-hand side vector), then taking the Chebyshev interval to be $\left[\frac{1.1M}{4},1.1M\right]$, where $M$ is the estimate of the largest eigenvalue of the relaxation-preconditioned linear system.  The factor of 1.1 (as used by default in PETSc) allows that $M$ is an under-estimate of the largest eigenvalue.  The factor of 4 chooses the Chebyshev polynomial to optimally damp the upper $3/4$ of the (estimated) spectrum of the relaxation-preconditioned linear system, as is suitable in 2D.  For relaxation on this problem, we use (additive) line relaxation.  Results at left of Table~\ref{tab:problem1} show results for line relaxation in both directions, while those at right present results for line relaxation only along lines parallel to the $x$-axis.

\begin{table}[!tbh]
  \begin{center}
  \scalebox{0.75}{
\begin{tabular}{|l|l|l|l|l|}
\hline
\multicolumn{5}{|c|}{Time To Solution (in minutes)}    \\ \hline
\backslashbox{$\varepsilon$}{N}  &   64   &     128   &     256   &     512   \\\hline
1.0e-04  &   9.70e-02 &   2.76e-01 &   9.85e-01 &   3.86e+00 \\
1.0e-06  &   1.12e-01 &   2.95e-01 &   1.01e+00 &   3.87e+00 \\
1.0e-08  &   1.14e-01 &   3.01e-01 &   1.01e+00 &   3.86e+00 \\
1.0e-10  &   1.18e-01 &   2.99e-01 &   1.01e+00 &   3.84e+00 \\
\hline
\multicolumn{5}{|c|}{Iteration counts}    \\ \hline
1.0e-04  &    3 &       3 &       4 &       4 \\
1.0e-06  &    3 &       3 &       4 &       4 \\
1.0e-08  &    3 &       3 &       3 &       4 \\
1.0e-10  &    3 &       3 &       3 &       3 \\
\hline
\multicolumn{5}{|c|}{Max Errors}    \\ \hline
1.0e-04  &    1.82e-02 &   8.10e-03 &   3.87e-03 &   1.91e-03 \\
1.0e-06  &    1.81e-02 &   8.00e-03 &   3.77e-03 &   1.85e-03 \\
1.0e-08  &    1.81e-02 &   7.98e-03 &   3.77e-03 &   1.85e-03 \\
1.0e-10  &    1.81e-02 &   7.98e-03 &   3.77e-03 &   1.85e-03 \\
\hline
\multicolumn{5}{|c|}{Energy Errors}    \\ \hline
1.0e-04  &    6.67e-02 &   3.89e-02 &   2.22e-02 &   1.25e-02 \\
1.0e-06  &    6.64e-02 &   3.87e-02 &   2.21e-02 &   1.24e-02 \\
1.0e-08  &    6.64e-02 &   3.87e-02 &   2.21e-02 &   1.24e-02 \\
1.0e-10  &    6.64e-02 &   3.87e-02 &   2.21e-02 &   1.24e-02 \\
\hline
\multicolumn{5}{|c|}{SD Errors}    \\ \hline
1.0e-04  &    6.68e-02 &   3.89e-02 &   2.22e-02 &   1.25e-02 \\
1.0e-06  &    6.65e-02 &   3.87e-02 &   2.21e-02 &   1.24e-02 \\
1.0e-08  &    6.65e-02 &   3.87e-02 &   2.21e-02 &   1.24e-02 \\
1.0e-10  &    6.65e-02 &   3.87e-02 &   2.21e-02 &   1.24e-02 \\
\hline
\end{tabular}

\begin{tabular}{|l|l|l|l|l|}
\hline
\multicolumn{5}{|c|}{Time To Solution (in minutes)}    \\ \hline
\backslashbox{$\varepsilon$}{N}  &   64   &     128   &     256   &     512   \\\hline
1.0e-04  &   7.18e-02 &   1.71e-01 &   5.73e-01 &   2.14e+00 \\
1.0e-06  &   8.42e-02 &   1.84e-01 &   5.76e-01 &   2.15e+00 \\
1.0e-08  &   8.72e-02 &   1.96e-01 &   5.74e-01 &   2.13e+00 \\
1.0e-10  &   8.79e-02 &   1.93e-01 &   5.75e-01 &   2.11e+00 \\
\hline
\multicolumn{5}{|c|}{Iteration counts}    \\ \hline
1.0e-04 &       2 &       3 &       3 &       4 \\
1.0e-06 &       2 &       2 &       3 &       4 \\
1.0e-08 &       2 &       2 &       2 &       3 \\
1.0e-10 &       2 &       2 &       2 &       2 \\
\hline
\multicolumn{5}{|c|}{Max Errors}    \\ \hline
1.0e-04  &   1.82e-02 &   8.10e-03 &   3.87e-03 &   1.91e-03 \\
1.0e-06  &   1.81e-02 &   7.97e-03 &   3.77e-03 &   1.85e-03 \\
1.0e-08  &   1.81e-02 &   7.97e-03 &   3.73e-03 &   1.85e-03 \\
1.0e-10  &   1.81e-02 &   7.97e-03 &   3.73e-03 &   1.81e-03 \\
\hline
\multicolumn{5}{|c|}{Energy Errors}    \\ \hline
1.0e-04 &   6.67e-02 &   3.89e-02 &   2.22e-02 &   1.25e-02 \\
1.0e-06 &   6.64e-02 &   3.87e-02 &   2.21e-02 &   1.24e-02 \\
1.0e-08 &   6.64e-02 &   3.87e-02 &   2.21e-02 &   1.24e-02 \\
1.0e-10 &   6.64e-02 &   3.87e-02 &   2.21e-02 &   1.24e-02 \\
\hline
\multicolumn{5}{|c|}{SD Errors}    \\ \hline
1.0e-04 &   6.68e-02 &   3.89e-02 &   2.22e-02 &   1.25e-02 \\
1.0e-06 &   6.65e-02 &   3.87e-02 &   2.21e-02 &   1.24e-02 \\
1.0e-08 &   6.65e-02 &   3.87e-02 &   2.21e-02 &   1.24e-02 \\
1.0e-10 &   6.65e-02 &   3.87e-02 &   2.21e-02 &   1.24e-02 \\
\hline
\end{tabular}
}
  \caption{Computational time, iteration counts, and errors for model problem in~\eqref{eq4}.  At left, results using line relaxation in both directions.  At right, results using line relaxation only along lines parallel to the $x$-axis.}\label{tab:problem1}
\end{center}
\end{table}

Computational times and iteration counts reported in Table~\ref{tab:problem1} show the hoped-for scaling, with a number of iterations that is bounded independently of both $\varepsilon$ and $N$, and a computational time that scales like $\mathcal{O}(N^2)$, the total number of grid-points in the mesh, for both relaxation strategies.  We note that the CPU times are somewhat larger than desirable for systems of these sizes.  This is due to a considerable time required for setting up the block relaxation scheme.  For a single problem on the $512\times 512$ grid using line relaxation in both directions, we find that over 95\% of the CPU time is taken by the preconditioner setup phase, including both initializing the multigrid hierarchy and computing the components of the block relaxation scheme.  About 95\% of this (90\% of the total CPU times) is spent in the PCPatch setup phase, for the block relaxation.  Of this, about 75\% is spent computing the index sets to specify the patches (over 65\% of the total CPU time).  Finally, almost 90\% of index set computation comes from retrieving the coordinates of the degrees of freedom from the PETSc DM (Domain Manager).  Thus, we see that almost 60\% of our total solve time is dedicated to a single line of code in Firedrake's interface to PCPatch.  In comparison, the FGMRES solve phase accounts for only about 2\% of the total runtime, about 5 seconds.  We note two important points about this.  First of all, much of the setup computation for the patches could be saved and reused, since the indexing of these patches relative to the global grids is consistent, independent of the number of levels in the multigrid hierarchy and the perturbation parameter, $\varepsilon$, from which the Shishkin meshes are computed.  Secondly, we have not systematically looked at optimizing this part of the computation, for example, by providing the DoF coordinates in a more accessible manner.

Comparing the left and right sides of Table~\ref{tab:problem1}, we see that both lead to scalable solvers for these problems, with low iteration counts to solution.  All three measures of the error in the final approximation are essentially identical between the two solvers.  From the analysis above, we expected line relaxation only along lines parallel to the $x$-axis to yield a scalable solver, and these results confirm this.  Moreover, we see the true benefit of this insight, since the computational cost of line relaxation is the most substantial cost of the algorithm.  Thus, using line relaxation only along lines in one direction leads to an algorithm that yields equally good approximate solutions in about half the computational time.

Finally, we note several important scalings in the error measures reported in Table~\ref{tab:problem1}.  Since we are computing on a Shishkin mesh, and measuring convergence in ``stronger'' norms than just the $L^2(\Omega)$ norm, we generally expect convergence like $\mathcal{O}(N^{-1}\ln N)$, leading to a reduction in error of a factor of about 0.56 between $N=256$ and $N=512$.  We see exactly such reductions in the Energy and SD Errors, and slightly better (consistent with $\mathcal{O}(N^{-1})$) reduction for the Max Error.  We also note that all reported errors are bounded independently of $\varepsilon$, as expected for a parameter-robust discretization, and that the SD Errors are nearly identical to the Energy Errors.  This is partly to be expected, since the weighting on the additional term in~\eqref{eq:SDError} is zero in the left-hand boundary layer and $\mathcal{O}(h_{\mathcal{T}})$ in the remainder of the domain; however, it confirms that there is not a significant error in streamline gradient term where the mesh spacing is coarsest.  We also note that the Energy and SD Errors reported in this table are consistent with those reported in Table I of Franz et al.~\cite{FRANZ20081818}.

Our second model problem is: 
\begin{equation}
-\varepsilon\Delta u - u_x+u=f ~~\text{in} ~~ \Omega=(0,1)^2,~~ u=0~~\text{on}~~ \partial\Omega \label{eq5}
\end{equation}
with exact solution chosen to be 
\begin{equation*}
 u(x,y)= \sin (\pi x) \dfrac{(1-e^{-y/\sqrt{\varepsilon}})(1-e^{-(1-y)/\sqrt{\varepsilon}})}{1-e^{-1/\sqrt{\varepsilon}}}
\end{equation*}
and a corresponding right-hand side function, $f$.  We note that this problem has only parabolic layers, so we consider a discretization on a uniform mesh in the $x$-direction and Shishkin mesh in the $y$-direction.  We make several adjustments in the discretization and solver setup to account for the lack of an exponential layer.  First, we now choose the SUPG stabilization parameter to be $h_{\mathcal{T}}$ for all cells not within a parabolic layer, and to be $h_{\mathcal{T}}^{4/3}$ for those cells in the layers near $y=0$ and $y=1$ (inspired by the discussion in Franz et al.~\cite{FRANZ20081818}).  Secondly, we slightly adapt the stopping tolerance for FGMRES, requiring the absolute residual norm to be reduced below $\sqrt{\varepsilon}/N^2$, as an $\varepsilon$-independent stopping tolerance did not quite yield robust results.  As above, we only use block relaxation along lines parallel to the $x$-axis.  The analysis in Section~\ref{ssec:identifying} suggests that, due to the lack of grid adaptation in the $x$-direction, line relaxation may be more than is required for robustness in this setting, but we were unable to identify suitable relaxation parameters to yield a scalable solver with point relaxation for this problem.  This likely reflects that the $\mathcal{O}(h)$ stencil entries from the convection term are dominant here, since the diffusive term contributes only stencil entries that are $\mathcal{O}(\varepsilon)$ and $\mathcal{O}(\sqrt{\varepsilon})$ in size.

Table~\ref{tab:problem2} presents solution times, solver iteration counts, and error measures for this problem.  Again, we see that solver time scales linearly with problem size, and that general solution times are smaller than in Table~\ref{tab:problem1}, due to the fact that we construct only a single set of blocks for the relaxation.  Iteration counts are again stable (but slightly higher than in Table~\ref{tab:problem1}, due to the stricter stopping tolerance).  Maximum norm errors are again independent of $\varepsilon$, but now are scaling consistently with $\mathcal{O}(N^{-2}\ln^2N)$, as might be expected since the relevant error estimates are dominated by the layer terms in the $x$-direction, which are not present in this example.  Similarly, we see that while the energy norm decreases like $\mathcal{O}(N^{-1}\ln N)$ for fixed $\varepsilon$, it now decays like $\varepsilon^{1/4}$ for fixed $N$.  The SD errors, in contrast, decay much more slowly in $\varepsilon$ for fixed $N$.  This leads to the interesting phenomenon that they decay like $\mathcal{O}(N^{-1}\ln N)$ for $\varepsilon = 10^{-4}$ (where they are quite similar to the energy-norm errors), but like $\mathcal{O}(N^{-2}\ln^2 N)$ for $\varepsilon = 10^{-10}$, where the streamwise gradient term is clearly the dominant source of error.  We note again that the energy-norm error measures presented here are consistent with those presented by Franz et al.~\cite{FRANZ20081818}.

\begin{table}[!tbh]
  \begin{center}
  \scalebox{0.75}{
\begin{tabular}{|l|l|l|l|l|}
\hline
\multicolumn{5}{|c|}{Time To Solution (in minutes)}    \\ \hline
\backslashbox{$\varepsilon$}{N}   &     64   &     128   &     256   &     512   \\\hline
1.0e-04  &    7.17e-02 &   1.77e-01 &   5.80e-01 &   2.21e+00 \\
1.0e-06  &    8.90e-02 &   1.92e-01 &   5.86e-01 &   2.21e+00 \\
1.0e-08  &    8.67e-02 &   1.97e-01 &   5.86e-01 &   2.21e+00 \\
1.0e-10  &    8.73e-02 &   1.93e-01 &   5.87e-01 &   2.17e+00 \\
\hline

\multicolumn{5}{|c|}{Iteration counts}    \\ \hline
1.0e-04  &       4 &       5 &       5 &       6 \\
1.0e-06  &       6 &       5 &       5 &       6 \\
1.0e-08  &       6 &       5 &       4 &       6 \\
1.0e-10  &       7 &       4 &       4 &       4 \\
\hline
\multicolumn{5}{|c|}{Max Errors}    \\ \hline
1.0e-04  &   1.07e-02 &   3.48e-03 &   1.07e-03 &   3.22e-04 \\
1.0e-06  &   1.08e-02 &   3.53e-03 &   1.10e-03 &   3.36e-04 \\
1.0e-08  &   1.08e-02 &   3.53e-03 &   1.11e-03 &   3.36e-04 \\
1.0e-10  &   1.08e-02 &   3.53e-03 &   1.11e-03 &   3.38e-04 \\
\hline
\multicolumn{5}{|c|}{Energy Errors}    \\ \hline
1.0e-04  &   1.32e-02 &   7.74e-03 &   4.42e-03 &   2.49e-03 \\
1.0e-06  &   4.19e-03 &   2.45e-03 &   1.40e-03 &   7.86e-04 \\
1.0e-08  &   1.33e-03 &   7.74e-04 &   4.42e-04 &   2.49e-04 \\
1.0e-10  &   4.24e-04 &   2.45e-04 &   1.40e-04 &   7.87e-05 \\
\hline
\multicolumn{5}{|c|}{SD Errors}    \\ \hline
1.0e-04  &   1.41e-02 &   7.91e-03 &   4.46e-03 &   2.49e-03 \\
1.0e-06  &   7.14e-03 &   3.18e-03 &   1.57e-03 &   8.26e-04 \\
1.0e-08  &   6.02e-03 &   2.22e-03 &   8.57e-04 &   3.59e-04 \\
1.0e-10  &   5.90e-03 &   2.09e-03 &   7.48e-04 &   2.72e-04 \\
\hline
\end{tabular}
  }
\caption{Computational time, iteration counts, and errors for model problem in~\eqref{eq5}}\label{tab:problem2}
\end{center}
\end{table}

\subsection{Semi-Structured Curvilinear Meshes}
\begin{figure}[t]
\begin{minipage}{.47\textwidth}
\centering
\includegraphics[trim={2cm 1cm 1cm 1cm},clip,scale=.52]{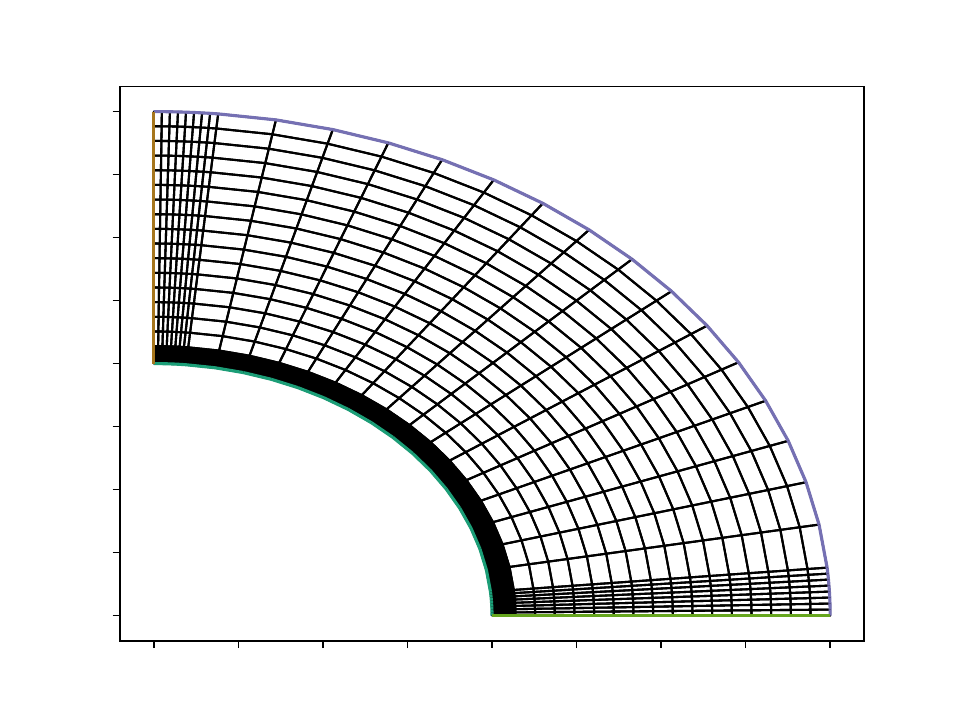}
\end{minipage}
\begin{minipage}{.5\textwidth}
\centering
\includegraphics[trim={4cm 10cm 1cm 10cm},clip,scale=.32]{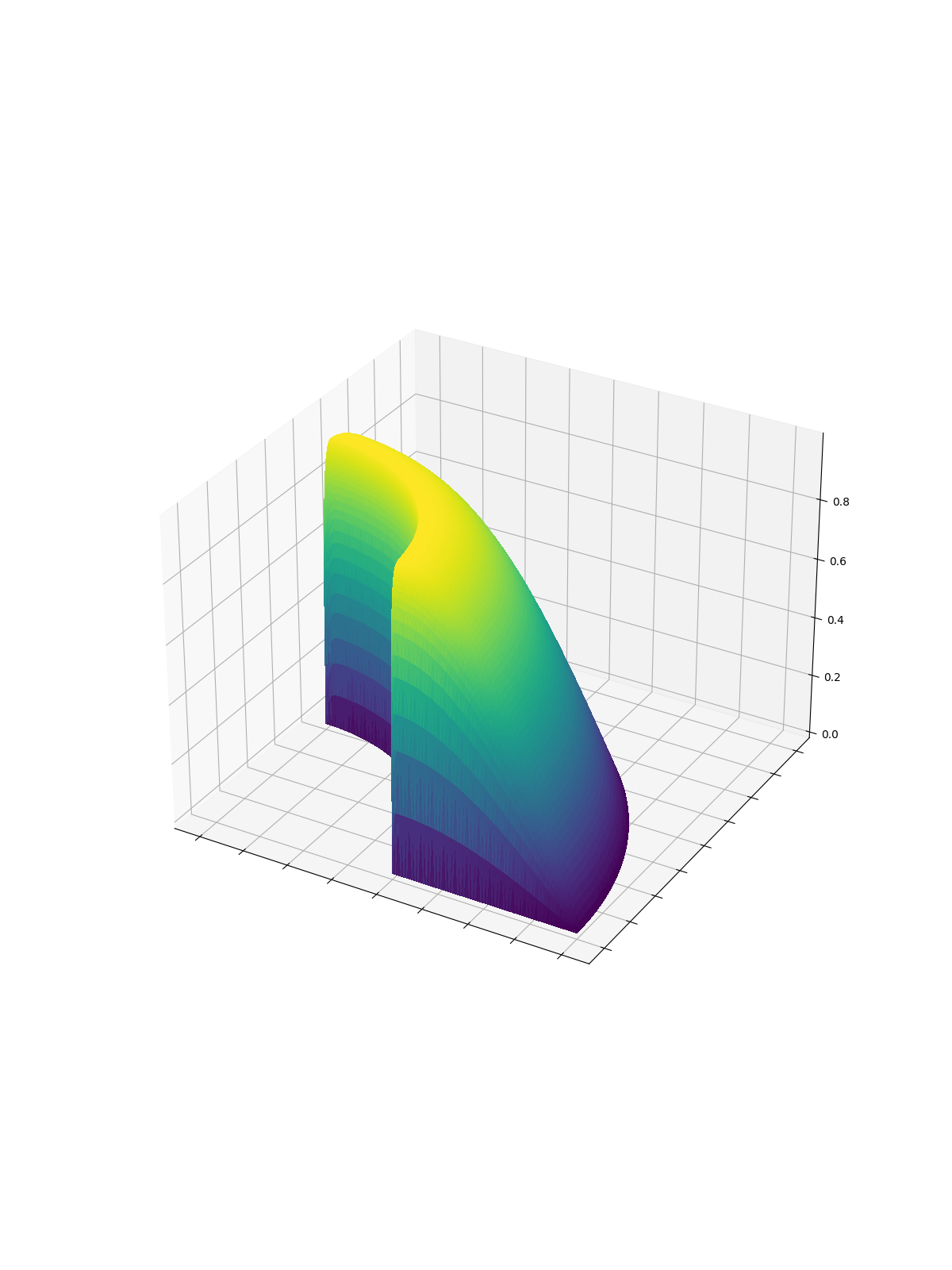}
\end{minipage}
\caption{The curvilinear mesh and solution at $\varepsilon=10^{-3}$}\label{fig:quart_ann}
\end{figure}

The main focus of our work on generalizing line relaxation in Firedrake was to be able to treat much more complicated model problems than could be solved using traditional multigrid line-relaxation codes.  In this section, we demonstrate that capability on a curved domain, taken to be a quarter of an annulus, with inner radius 1 and outer radius 2.  We realize such domains using a mesh-mapping approach, where we first construct a hierarchy of nested layer-adapted meshes of the unit square, then map each mesh in the hierarchy to the quarter-annulus domain, with polar coordinates $r = x+1$, $\theta = y\pi/2$, so that the adaptation along $x=0$ becomes adaptation along the inner radius of the annulus, while that along $y=0$ and $1$ becomes adaptation along the edges $\theta = 0$ and $\theta = \pi/2$ of the quarter-annulus domain, respectively.  On this domain, we consider the PDE
\begin{equation}
-\varepsilon \Delta u - (r+1)\frac{\partial u}{\partial r} + u = f,\label{eq6}
\end{equation}
with zero Dirichlet boundary conditions, so that the convection term is tangent to the radial vector and pointing towards the origin, leading to the expectation of parabolic layers along $\theta = 0$ and $\theta = \pi/2$ and an exponential layer at the inner radius of $r=1$.  An important aspect of achieving accuracy on curved domains is adequately resolving the domain curvature in the mesh.  Here, we use a second-order piecewise polynomial basis for the domain geometry, which is sufficient to ensure that our mesh perfectly resolves the geometry of the domain.  An example mesh with 32 elements in each direction (radially and azimuthally) for $\varepsilon = 10^{-3}$ is shown at left of Figure~\ref{fig:quart_ann}.

To test the effectiveness of our solvers on this domain, we consider a similar mapping of the model problem from the first unit-square example considered above to polar coordinates on this domain, with a manufactured solution given by
\begin{equation}
u(r,\theta)=
\left(\cos\left(\frac{\pi(r-1)}{2}\right) - \frac{e^{-(r-1)/\varepsilon}- e^{-1/\varepsilon}}{1-e^{-1/\varepsilon}}\right) \frac{(1-e^{-(2\theta/\pi)/\sqrt{\varepsilon}})(1-e^{-(1-2\theta/\pi)/\sqrt{\varepsilon}})}{1-e^{-1/\sqrt{\varepsilon}}}.
\label{eq_crv}
\end{equation}
A computed solution for the resulting discrete problem is shown at right of Figure~\ref{fig:quart_ann}.

One key question in this work is determining the transition points in the original Shishkin mesh on the unit square to ensure suitable finite-element convergence for such solutions on the mapped mesh.  There is no significant distortion in lengths for the mapping $r = x-1$, so we keep the standard Shishkin transition point, $\lambda_1$, defined in~\eqref{eq:xShishkin_transition}.  Lengths in the $\theta$ direction, however, are notably distorted by the mapping from $y$ to $\theta = y\pi/2$.  Through experimentation, we found the best results were obtained by scaling the second term in $\lambda_2$ from~\eqref{eq:yShishkin_transition} by a factor of $2.75/\pi$, slightly larger than the naive scaling of $2/\pi$ that we tried first.  Here, taking representative ``interior'' mesh size, $h$, to be the square-root of the volume of the largest element in the mesh, we use SUPG stabilization parameter of $0.55h$ in the interior region, and 0 in the layer regions.

Very few adaptations to the multigrid preconditioning framework developed above are required for these examples.  We continue to use natural finite-element interpolation operators, now on the mapped meshes, and use weighted Chebyshev acceleration of relaxation as described above.  As discussed above, successful solvers can be developed using only line relaxation along lines of fixed $\theta$ value, since these are the strongest connections in the stencil.  For our first experiments, shown at left of Table~\ref{tab:problem_curv}, we define the divisions for PCPatch so that each line of nodes with common $\theta$-value in the mesh makes one patch.  These results are computed for a piecewise linear discretization with V(2,2) cycles and an FGMRES stopping criterion requiring a reduction in the residual norm below an absolute factor of $\varepsilon/N$.

\begin{table}[!tbh]
  \begin{center}
  \scalebox{0.75}{
  \begin{tabular}{c c}
 \begin{tabular}{|l|l|l|l|l|}
\hline
\multicolumn{5}{|c|}{Time To Solution (in minutes)}    \\ \hline
\backslashbox{$\varepsilon$}{N}   &     64   &     128   &     256   &     512   \\\hline
1.0e-02  &    2.75e-01 &   8.07e-01 &   3.04e+00 &   1.26e+01 \\
1.0e-03  &    2.65e-01 &   8.74e-01 &   3.22e+00 &   1.33e+01 \\
1.0e-04  &    2.77e-01 &   9.11e-01 &   4.04e+00 &   1.37e+01 \\
1.0e-05  &    2.71e-01 &   9.12e-01 &   3.71e+00 &   1.43e+01 \\
1.0e-06  &    3.07e-01 &   9.73e-01 &   3.85e+00 &   1.67e+01 \\
\hline

\multicolumn{5}{|c|}{Iteration counts}    \\ \hline
1.0e-02  &        3 &       3 &       3 &       4 \\
1.0e-03  &        2 &       3 &       3 &       4 \\
1.0e-04  &        3 &       4 &       6 &       5 \\
1.0e-05  &        3 &       4 &       5 &       6 \\
1.0e-06  &        4 &       5 &       6 &       8 \\
\hline
\multicolumn{5}{|c|}{Max Errors}    \\ \hline
1.0e-02  &    3.03e-03 &   1.10e-03 &   5.37e-04 &   1.57e-04 \\
1.0e-03  &    1.04e-02 &   2.84e-03 &   7.14e-04 &   1.74e-04 \\
1.0e-04  &    1.05e-02 &   3.30e-03 &   1.06e-03 &   3.32e-04 \\
1.0e-05  &    1.05e-02 &   3.30e-03 &   1.07e-03 &   3.33e-04 \\
1.0e-06  &    1.05e-02 &   3.30e-03 &   1.07e-03 &   3.34e-04 \\
\hline
\multicolumn{5}{|c|}{Energy Errors}    \\ \hline
1.0e-02  &    5.63e-02 &   3.28e-02 &   1.87e-02 &   1.05e-02 \\
1.0e-03  &    6.45e-02 &   3.75e-02 &   2.13e-02 &   1.20e-02 \\
1.0e-04  &    6.59e-02 &   3.84e-02 &   2.19e-02 &   1.23e-02 \\
1.0e-05  &    6.63e-02 &   3.87e-02 &   2.21e-02 &   1.24e-02 \\
1.0e-06  &    6.68e-02 &   3.88e-02 &   2.21e-02 &   1.25e-02 \\
\hline
\multicolumn{5}{|c|}{SD Errors}    \\ \hline
1.0e-02  &    5.64e-02 &   3.28e-02 &   1.87e-02 &   1.05e-02 \\
1.0e-03  &    6.45e-02 &   3.75e-02 &   2.13e-02 &   1.20e-02 \\
1.0e-04  &    6.61e-02 &   3.85e-02 &   2.20e-02 &   1.23e-02 \\
1.0e-05  &    6.66e-02 &   3.87e-02 &   2.21e-02 &   1.24e-02 \\
1.0e-06  &    6.71e-02 &   3.88e-02 &   2.22e-02 &   1.25e-02 \\
\hline
\end{tabular}
&

\begin{tabular}{|l|l|l|l|l|}
\hline
\multicolumn{5}{|c|}{Time To Solution (in minutes)}    \\ \hline
\backslashbox{$\varepsilon$}{N}   &     64   &     128   &     256   &     512   \\\hline
1.0e-02  &    3.79e-01 &   1.36e+00 &   6.07e+00 &   2.91e+01 \\
1.0e-03  &    4.07e-01 &   1.40e+00 &   4.86e+00 &   2.32e+01 \\
1.0e-04  &    4.27e-01 &   1.40e+00 &   4.93e+00 &   2.10e+01 \\
1.0e-05  &    3.93e-01 &   1.33e+00 &   4.74e+00 &   1.98e+01 \\
1.0e-06  &    4.45e-01 &   1.30e+00 &   4.51e+00 &   1.94e+01 \\
\hline
\multicolumn{5}{|c|}{Iteration counts}    \\ \hline
1.0e-02  &        4 &       5 &       7 &      10 \\
1.0e-03  &        4 &       5 &       5 &       7 \\
1.0e-04  &        4 &       5 &       5 &       6 \\
1.0e-05  &        4 &       5 &       5 &       6 \\
1.0e-06  &        6 &       5 &       5 &       6 \\
\hline
\multicolumn{5}{|c|}{Max Errors}    \\ \hline
1.0e-02  &    5.21e-05 &   5.10e-04 &   7.45e-07 &   1.44e-04 \\
1.0e-03  &    1.10e-04 &   7.35e-04 &   1.21e-05 &   2.69e-04 \\
1.0e-04  &    2.29e-04 &   2.69e-05 &   9.90e-06 &   1.34e-06 \\
1.0e-05  &    2.57e-04 &   6.62e-04 &   3.46e-05 &   8.60e-06 \\
1.0e-06  &    2.39e-03 &   1.51e-04 &   3.07e-05 &   7.49e-06 \\
\hline
\multicolumn{5}{|c|}{Energy Errors}    \\ \hline
1.0e-02  &    4.20e-03 &   1.54e-03 &   4.68e-04 &   3.31e-04 \\
1.0e-03  &    4.83e-03 &   1.70e-03 &   5.33e-04 &   2.61e-04 \\
1.0e-04  &    5.00e-03 &   1.71e-03 &   5.58e-04 &   1.77e-04 \\
1.0e-05  &    5.04e-03 &   1.71e-03 &   5.58e-04 &   1.76e-04 \\
1.0e-06  &    9.19e-03 &   1.77e-03 &   5.58e-04 &   1.76e-04 \\
\hline
\multicolumn{5}{|c|}{SD Errors}    \\ \hline
1.0e-02  &    4.20e-03 &   1.65e-03 &   4.68e-04 &   3.63e-04 \\
1.0e-03  &    4.83e-03 &   2.03e-03 &   5.33e-04 &   3.45e-04 \\
1.0e-04  &    5.00e-03 &   1.71e-03 &   5.58e-04 &   1.77e-04 \\
1.0e-05  &    5.15e-03 &   1.83e-03 &   5.60e-04 &   1.78e-04 \\
1.0e-06  &    9.19e-03 &   1.77e-03 &   5.60e-04 &   1.78e-04 \\
\hline
\end{tabular} 
  \end{tabular}
  }
\caption{Computational time, iteration counts, and errors for model problem in~\eqref{eq6}.  At left, results for piecewise linear finite-element solutions.  At right, results for piecewise quadratic finite-element solutions.}\label{tab:problem_curv}
\end{center}
\end{table}

Considering the results at left of Table~\ref{tab:problem_curv}, we see quite comparable results to those in Table~\ref{tab:problem1} for a rectangular-grid version of this problem.  The iteration counts reported here are slightly higher than those in Table~\ref{tab:problem1}, due to the stricter stopping criterion needed to achieve finite-element convergence in this case.  The overall computing times, however, are roughly a factor of four larger.  This is due to the use of a piecewise quadratic mesh representation, which increases the dominant cost of assembling the lines for relaxation based on the underlying mesh data structure, for which optimization has not yet been considered.  Errors in the discrete maximum norm scale like $\mathcal{O}(N^{-2}\ln^2(N))$, slightly better than observed in the unit square case, while those in the energy and SD norms scale like $\mathcal{O}(N^{-1}\ln(N))$.  As above, the SD errors are nearly identical to the Energy errors, confirming there is not a significant error in the streamline gradient term where the mesh spacing is coarse.

The mesh mapping approach used here also allows us to explore the impacts on accuracy of using higher-order discretizations, such as the isoparametric discretization that arises when using piecewise quadratic finite-elements with the piecewise-quadratic mesh representation.  Results for this are shown at right of Table~\ref{tab:problem_curv}.  Here, we slightly adjust the FGMRES stopping tolerance, now requiring an absolute residual norm below $0.1\varepsilon/N^2$ and using V(3,3) multigrid cycles as a preconditioner.  We adapt our line relaxation definition so that the number of patches stays equal to the number of nodes of our mesh in the $\theta$-direction, grouping degrees of freedom that share the same $\theta$ value at their location (at a node, edge, or cell midpoint of a mesh), and pairing these so that each line now includes the nodes (as used in piecewise linear case), DoFs at edges between those nodes, and those edge and cell-centered DoFs on the adjacent line of fixed $\theta$ with $\theta$-value smaller than those at the nodes.  These lines, consequently have four times as many DoFs as those in the piecewise linear case.  The results at right of Table~\ref{tab:problem_curv} show again that the setup costs are dominant, as we have problems with four times as many degrees of freedom, but only computing times increasing by a factor of 2-3 (or less).  Corresponding increases in accuracy are quite clear, with observed convergence of $\mathcal{O}(N^{-2})$ or better in the discrete maximum norm and $\mathcal{O}(N^{-2}\ln^2(N))$ in the energy and SD norms.  We note the SD norms are still generally dominated by the energy error, although this is slightly less clear for large $\varepsilon$.

\subsection{The Hemker Problem}\label{ssec:hemker}

\begin{figure}
\centering
\includegraphics[width=0.8\textwidth]{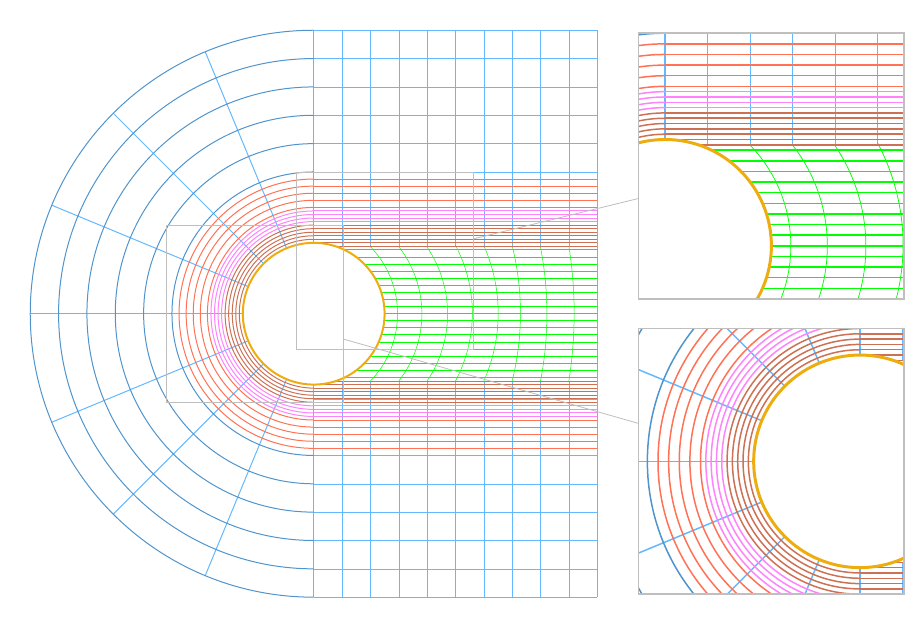}
\caption{At left, sketch of the domain, $\Omega$, and its mesh for the solution of the Hemker problem.  At right, a zoom-in of the structure of the three nested boundary layers near the bottom of the unit circle.}\label{fig:Hemker}
\end{figure}

The \emph{Hemker problem}~\cite{hemker1996singularly} was proposed by Piet Hemker as ``a relatively simple problem on a non-trivial domain, so that the problems related with a proper meshing and with the treatment of the boundary conditions can be studied more thoroughly.''  The problem proposed is a simple two-dimensional convection-diffusion problem; however, it is posed on the exterior of the unit circle.  While Hemker considered the problem on an infinite domain, where a series solution and its asymptotic expansion for small $\varepsilon$ can be found, we will follow the more recent presentation by Hegarty and O'Riordan~\cite{hegarty2020numerical} on a finite domain.  This problem has been considered in several variations in recent years~\cite{augustin2011assessment, john2013robust, barrenechea2018unified}.
Here, we consider the continuum domain
\begin{equation*}
\Omega = \{(x,y) \mid 1 \leq x^2+y^2 \leq 4^2\}\cup \{(x,y) \mid -4 \leq y \leq 4, 0 < x < 4, 1 < x^2+y^2\},
\end{equation*}
pictured at left of Figure~\ref{fig:Hemker}, and the boundary value problem
\begin{align*}
-\varepsilon \Delta u +  u_x &= 0 \text{ in }\Omega \\
 u(x,y)&=1,~\text{if} ~ x^2+y^2=1~~~~ \\
 u(x,y)&=0,~\text{if} ~ x^2+y^2=4~\text{and}~x\leq0~~ \\ 
 u(x,y)&=0,~\text{if} ~ y=\pm4~\text{and}~0\leq x\leq 4~~ \\ 
 u_x(x,y)&=0,~\text{if} ~ -4< y<4~\text{and}~x=4.~~ 
\end{align*}

The nonzero boundary condition on the inner circle, coupled with the singularly perturbed nature of the governing PDE lead to a complicated layer structure in $u(x,y)$ that has been well-studied.  Theorem 4 of Hegarty and O'Riordan~\cite{hegarty2020numerical} establishes three bounds on the solution decay that inform the layer structure of the solution and, in turn, lead to the structure of the Shishkin-style mesh that we adopt.  First of all, in the left region of the domain, bounded between the two half-circles with radii 1 and 4, the solution decays like $e^{cos(\theta)(r-1)/\varepsilon}$ (where $(r,\theta)$ are standard polar coordinates in this region), noting that $\cos(\theta) \leq 0$ in this region.  This leads to the need to form a radial mesh with spacing that is $\mathcal{O}(\varepsilon)$ near the $r=1$ boundary of this region.  Secondly, when $|y|\geq 1$, the solution is bounded by a constant times $e^{-(|y|-1)/\sqrt{\varepsilon}}$.  In the right half plane, this leads to the need to form a mesh with spacing that is $\mathcal{O}(\sqrt{\varepsilon})$ near the lines $y = \pm 1$, both near the boundary where $x^2+y^2=1$ and in the interior of the domain.  Finally, near $x=0$, there is a layer of width $\varepsilon^{1/3}$ in $x$ and $\varepsilon^{2/3}$ in $y$.  Again, this requires meshing with spacing proportional to $\varepsilon^{2/3}$ in $y$ near the boundary of the inner unit circle.  We note that the convection direction in this example is towards the right (in contrast to the earlier examples and discussion), but that the presence of a Neumann boundary condition on the right edge of the domain implies that no boundary layers arise there.

While Hegarty and O'Riordan make use of an overlapping pair of discretization meshes, to separately resolve the boundary layer in the left part of the domain and the interior layer in the right~\cite{hegarty2020numerical}, we consider a single mesh, with additional grading to attempt to resolve all of the layer structure described above with minimal \textit{over-resolution}, where small elements are included with little computational benefit.  In the left half of the domain, we consider a mesh with grid points chosen based on a tensor-product mesh in polar coordinates, with even spacing in the angular direction, with $N$ angular intervals for $\pi/2 \leq \theta \leq 3\pi/2$.  In this part of the domain, in the radial direction, we define three transition points, $\sigma_1$, $\sigma_2$, and $\sigma_3$, in order to divide the radial interval $1 \leq r \leq 4$ into four regions, $[1,1+\sigma_1]$, $[1+\sigma_1,1+\sigma_2]$, $[1+\sigma_2,1+\sigma_3]$, $[1+\sigma_3,4]$, with $N/4$ radial intervals in the first three regions, and $N/2$ intervals in the last region.  Here, we take
\begin{equation*}
\sigma_1:=\min\left\lbrace 0.25, \varepsilon \ln N\right\rbrace~~
\sigma_2:=\min\left\lbrace 0.3,\varepsilon^{2/3} \ln N\right\rbrace~~
\sigma_3:=\min\left\lbrace 0.35,\varepsilon^{1/2} \ln N\right\rbrace~~
\end{equation*}
In the right half of the domain, we have $N/2$ evenly spaced mesh intervals in the horizontal direction for $0\leq x \leq 4$ along both the top and bottom edges of the domain, and these are ``wrapped around'' the unit circle for $-1 \leq y \leq 1$.  In the vertical direction, mesh intervals are continued from the radial divisions in the left half-plane, leading to layers structures for $|y|\geq 1$ that match the transition points and interval distributions described above.  For $-1\leq y \leq 1$, we use $3N/2$ intervals in the vertical direction.  We distribute these elements in 3 groups, drawing horizontal lines connecting the unit circle to the right-hand boundary at $y = \pm (1-\sigma_3)$, and using $N/2$ elements in the vertical direction in each of the three subregions thus created.  We note that the resulting elements do not have edges perfectly aligned with the $x$-axis, as gmsh naturally divides the elements by arc length on the unit circle and uniformly-in-$y$ on the right-hand boundary, leading to slightly ``tilted'' edges in the horizontal direction.  We note that while the solution is not expected to exhibit significant layer structure in the region to the right of the unit circle, we need some mesh resolution here to simply resolve the geometry of the circle itself.  A hierarchy of such meshes is constructed by using gmsh~\cite{geuzaine2009gmsh} to create a coarse mesh that resolves the transition points and the geometry, then using uniform (Cartesian) refinement of that mesh, then moving the resulting mesh points to properly resolve the polar geometry as required.  This results in an almost nested hierarchy of grids, with improving resolution of the curved geometry as the grids are refined.  For geometric multigrid on this hierarchy, we simply extrapolate in cases where fine-grid points are outside of the coarse-grid domain.  Due to a current software limitation, we are not able to construct a piecewise-quadratic mesh geometry starting from gmsh meshes, which will limit accuracy in the results to follow.

To discretize the problem, we use the SUPG discretization in~\eqref{eq:SUPG_weak}, with stabilization parameter chosen to be 0 for all elements contained in the polar-coordinates region $1\leq r \leq 1+\sigma_3$, and $0.55h_{\mathcal{T}}$ for all other elements, corresponding to the choices of $\tau_0 = 0.55$ and $\tau_1 = 0$ in~\eqref{Tau_0}.  We use multigrid-preconditioned FGMRES with a stopping tolerance of $1/N^2$, noting that this was found experimentally to produce errors that are generally within a few percent of those observed when using a direct solver, and no worse than a factor of 10\% larger than these.  The line relaxation for this problem is constructed using radial lines in the left half plane and vertical lines in the right half plane, to resolve the mesh refinement.  Here, we found that using the Chebyshev relaxation described above was less effective, so we use right-preconditioned GMRES as relaxation, with 3 steps for both pre- and post-relaxation.

Since we do not have an analytical solution to the Hemker problem on a finite domain, we use refined computations to compute reference solutions for comparisons.  While Hegarty and O'Riordan use the classical ``double-mesh method'' to estimate convergence~\cite{hegarty2020numerical}, comparing the solution on each mesh with fixed $N$ to that with mesh parameter $2N$, we use the natural finite-element analogue (as in Hill and Madden~\cite{NMTMA-14-559}) and compare the piecewise-linear solution on a given mesh with the piecewise-quadratic finite-element solution generated on the same mesh.  We compute the higher-order reference solution using a direct solver (MUMPS~\cite{amestoy2000mumps}), to avoid introducing bias into the error computations.  We note, however, that this reinforces meshing effects in the error calculations, since a high-order reference solution on a mesh that poorly resolves the domain cannot accurately represent the true PDE solution on the true domain.  For this domain, we note that there is an inherent difficulty in meshing along both the inner and outer radii of the annulus with piecewise linear geometry.  On the inner radius, $r=1$, any given mesh includes some area inside the unit circle, connecting two mesh points on the circle by a straight line.  On the outer radius, $r=4$, the opposite occurs, and all meshes with straight edges exclude part of the physical domain by connecting two points on the outer circle with a straight line.  We note that the effects of this, particularly in the numerical resolution of boundary layers near $r=1$, precludes this discretization achieving high precision as $\varepsilon \rightarrow 0$.  Preliminary experiments (using direct solvers) show that meshing errors become dominant for $\varepsilon \leq 10^{-3}$, so we exclude such results from this paper.  In future work, we will explore the use of isoparametric elements to better resolve the geometry near these curved boundaries. Another option is to consider other stabilizations of the piecewise linear finite-element discretization, such as those proposed in~\cite{John_etal_2022}, to yield better-quality approximations on the same meshes.

\begin{table}[!tbh]
  \begin{center}
  \scalebox{0.75}{
\begin{tabular}{|l|l|l|l|l|l|}
  \hline
\multicolumn{6}{|c|}{Time To Solve (in minutes)}    \\ \hline
\backslashbox{$\varepsilon$}{N}  & 32   & 64  & 128        & 256        & 512 \\\hline
1.0e+00  &    7.16e-02 &   1.74e-01 &   5.47e-01 &   2.06e+00 &    8.24e+00 \\
1.0e-01  &    7.44e-02 &   1.78e-01 &   5.50e-01 &   2.01e+00 &    8.04e+00 \\
1.0e-02  &    1.08e-01 &   1.91e-01 &   5.98e-01 &   2.17e+00 &    8.40e+00 \\
1.0e-03  &    7.18e-02 &   2.22e-01 &   8.54e-01 &  3.62e+00  &    1.54e+01  \\

\hline
\multicolumn{6}{|c|}{Degrees of Freedom}    \\ \hline
   & 3464 & 13584& 53792& 214080 &854144\\\hline
\hline
\multicolumn{6}{|c|}{Iteration counts}      \\ \hline
1.0e+00  &        5 &       7 &       9 &      13 &       17 \\
1.0e-01  &        4 &       5 &       7 &      10 &       14 \\
1.0e-02  &        8 &      12 &      14 &      18 &       21 \\
1.0e-03  &        3 &       7 &      18 &      33 &       40 \\

\hline
\multicolumn{6}{|c|}{Max Errors}      \\ \hline
1.0e+00  &    1.30e-02 &   7.71e-03 &   4.21e-03 &   2.19e-03 &    1.12e-03 \\
1.0e-01  &    2.70e-02 &   1.20e-02 &   5.86e-03 &   2.91e-03 &    1.44e-03 \\
1.0e-02  &    3.18e-02 &   1.28e-02 &   5.68e-03 &   2.82e-03 &    1.58e-03 \\
1.0e-03  &    7.01e-01 &   2.58e-01 &   7.01e-02 &   1.81e-02 &    7.21e-03 \\
\hline
	 \multicolumn{6}{|c|}{Energy Errors}      \\ \hline
1.0e+00  &    1.06e-01 &   5.91e-02 &   3.20e-02 &   1.72e-02 &    9.19e-03 \\
1.0e-01  &    1.12e-01 &   5.84e-02 &   3.04e-02 &   1.59e-02 &    8.30e-03 \\
1.0e-02  &    1.14e-01 &   6.72e-02 &   3.88e-02 &   2.20e-02 &    1.23e-02 \\
1.0e-03  &   5.63e-01  &   1.80e-01 &   5.20e-02 &   1.98e-02 &    7.30e-03 \\
\hline
	 \multicolumn{6}{|c|}{SD Errors}      \\ \hline
1.0e+00  &    1.09e-01 &   5.99e-02 &   3.22e-02 &   1.72e-02 &    9.20e-03 \\
1.0e-01  &    1.14e-01 &   5.89e-02 &   3.06e-02 &   1.59e-02 &    8.32e-03 \\
1.0e-02  &    1.14e-01 &   6.72e-02 &   3.88e-02 &   2.20e-02 &    1.23e-02 \\
1.0e-03  &    5.65e-01 &   1.80e-01 &   5.20e-02 &   1.98e-02 &    7.30e-03 \\
\hline
\end{tabular}
}
\caption{Computational time, iteration counts, and errors for Hemker Problem.}
\label{tab2}
\end{center}
\end{table}

Numerical results for the Hemker problem are presented in Table~\ref{tab2}.  We note that the computational times are somewhat increased from the unit-square domain problems, but this is because the meshing of the Hemker domain leads to larger meshes for the same value of $N$ (as reported in the table).  We note that results look very promising for both the discretization and solver for $\varepsilon = 1, 0.1$, and $0.01$, before degrading (as noted above) at $\varepsilon = 10^{-3}$.  Comparing to Table 5 of Hegarty and O'Riordan~\cite{hegarty2020numerical}, we see that the maximum-norm errors reported herein are improved in comparison to results reported there.  We note, however, that the different choices of mesh mean that direct comparison is not feasible.

\section{Conclusions} \label{sec:conclusions}
Despite a long history of research, many open questions remain in the efficient and accurate solution of singularly perturbed differential equations.  In this work, we continue previous research by Ramage~\cite{MR1715555} and Gaspar, Clavero, and Lisbona~\cite{gaspar2002some} into the development of efficient geometric multigrid solvers for singularly perturbed convection diffusion equations.  We show that these techniques can be applied very effectively to standard model problems on the unit square.  Moreover, we demonstrate that it is possible to extend the line relaxation paradigm from tensor-product meshes of the unit square to structured meshes on more complex domains, as demanded by the Hemker problem~\cite{hemker1996singularly} considered in Section~\ref{ssec:hemker}.

Two important issues remain for future work.  First, the block relaxation used here is implemented using PETSc and Firedrake, via the PCPatch interface.  This results in substantial runtimes due to slow data access for mesh coordinates through this interface.  In future work, we will consider adding information to the interface in order to substantially speed up this process.  Secondly, while the numerical results for the Hemker problem are promising, they clearly suffer at small $\varepsilon$, due to the geometric error in resolving the curved domain boundaries using first-order geometry (elements with straight edges).  In future work, we will extend the mesh-mapping approach used for the quarter annulus example to more complicated domains, such as this one.  With these improvements, we believe the resulting solvers could be extended to interesting and challenging models of two- and three-dimensional fluid flow in complex geometries, where boundary layers are known to play an important role.

\bibliography{pap1}

\end{document}